\documentclass[reqno,12pt]{amsart}

\usepackage{amscd}
\usepackage{amsfonts}
\usepackage{amsmath}
\usepackage{amssymb}
\usepackage[american, english]{babel}
\usepackage{bbm}
\usepackage{bookmark}
\usepackage{cmap}
\usepackage{dsfont}
\usepackage{enumerate}
\usepackage{epigraph}
\usepackage{euscript}
\usepackage[myheadings]{fullpage}
\usepackage{graphicx}
\usepackage{mathabx}
\usepackage{mathrsfs}
\usepackage{mathtools}
\usepackage{refcount}
\usepackage{stmaryrd}
\usepackage{thinsp}
\usepackage[safe]{tipa}
\usepackage{verbatim}
\usepackage{xr-hyper}
\usepackage{hyperref}
\usepackage[all]{xy}

\numberwithin{equation}{subsection}

\newtheorem{thm}{Theorem}[subsection]
\newtheorem*{thm*}{Theorem}
\newtheorem{cor}[thm]{Corollary}
\newtheorem{lmm}[thm]{Lemma}

\newtheorem{prop}[thm]{Proposition}
\newtheorem{prop-const}[thm]{Proposition-Construction}

\newtheorem{rem}[thm]{Remark}

\newtheorem{defi}[thm]{Definition}

\newcommand{\on}{\operatorname}
\newcommand{\Spec}{\on{Spec}}
\newcommand{\coker}{\on{coker}}

\title{On the $p$-supports of a holonomic $\mathcal{D}$-module}
\author[Thomas Bitoun]{Thomas Bitoun}\thanks{Mathematical Institute, University of Oxford, Oxford OX2 6GG, UK; tbitoun@gmail.com}
\date{}

\begin{document}

\maketitle

\begin{abstract} For a smooth variety $Y$ over a perfect field of positive characteristic, the sheaf $D_Y$ of crystalline differential operators on $Y$ (also called the sheaf of $PD$-differential operators) is known to be an Azumaya algebra over $T^*_{Y'},$ the cotangent space of the Frobenius twist $Y'$ of $Y.$  
Thus to a sheaf of modules $M$ over $D_Y$ one can assign a closed subvariety of $T^*_{Y'},$ called the $p$-support, namely the support of $M$ seen as a sheaf on $T^*_{Y'}.$ We study here the family of $p$-supports assigned to the reductions modulo primes $p$ of a holonomic $\mathcal{D}$-module. We prove that the Azumaya algebra of differential operators splits on the regular locus of the $p$-support and that the $p$-support is a Lagrangian subvariety of the cotangent space, for $p$ large enough. The latter was conjectured by Kontsevich. Our approach also provides a new proof of the involutivity of the singular support of a holonomic $\mathcal{D}$-module, by reduction modulo $p.$

\end{abstract}

\setcounter{tocdepth}{1}

\section{Introduction}

Let $Y$ be a smooth variety over a perfect field. 
We may consider two  
sheaves of differential operators on $Y:$ on the one hand the sheaf $D^{(\infty)}_Y$ constructed by Grothendieck in EGA IV and on the other the sheaf $D_Y$ of crystalline differential operators, also called the sheaf of $PD$-differential operators, see e.g. \cite{BMR} and \cite{BO}. These sheaves coincide if the base field is of characteristic zero but they are very different if it is of positive characteristic. For example, $D^{(\infty)}_Y$ is typically a centerless sheaf of non-Noetherian rings. But $D_Y$ has a large center canonically identified with the symmetric algebra of the tangent sheaf of the Frobenius twist $Y'$ of $Y.$ And it is an Azumaya algebra over it, by \cite{BMR}. 

Thus to a 
coherent sheaf of modules $M$ over $D_Y$ one can assign a closed subvariety of the cotangent space $T^*_{Y'}$ called the $p$-support, see Definition \ref{defi: p-supp}. Namely it is the support of $M$ seen as a sheaf on $T^*_{Y'}.$ One might hope that the $p$-support is analogous to the classical notion of singular support of a say complex $\mathcal{D}$-module. The results of this paper confirm this hope. Note however that without further restrictions on $M,$ the $p$-support is  
arbitrary. Indeed every closed subvariety $Z$ of $T^*_{Y'}$ is the $p$-support of the quotient of $D_Y$ by the corresponding ideal $I_Z$ of its center. 

We consider here  
the reductions of a holonomic $\mathcal{D}$-module modulo large primes. 
That is to say, let $S$ be an integral scheme dominant and of finite type over $\Spec{\mathbb{Z}},$ for example the spectrum of the ring of integers of a number field, let $X$ be a smooth $S$-scheme of relative dimension $n.$ Then for $\mu$ the generic point of $S$ and for any closed point $s$ of $S,$ the generic fiber $X_\mu$ is a smooth variety over a field of characteristic zero while the fiber $X_s$ is a smooth variety over a field of positive characteristic. Let $D_{X/S}$ be the sheaf of relative crystalline differential operators and let $M$ be a coherent left $D_{X/S}$-module. Suppose that the generic fiber $M_\mu$ of $M$ is a nonzero holonomic $D_{X_\mu}$-module. For every closed point $s$ of $S,$ we let $M_s$ be the fiber of $M$ over $s.$ Our main results are:

\begin{enumerate}[(a)]
\item \label{item: purity} (Theorem \ref{thm:purity}) There is an open dense subset $U\subset S$ such that for all closed points $s$ of $U,$ the $p$-suppport of $M_s$ is equidimensional of dimension $n.$

\item \label{item: bound} (Theorem \ref{thm:bound}) Suppose that $X= \mathbb{A}^n_S$ with coordinates $\{x_1, \dots, x_n\}.$ For all closed points $s$ of $S$ of characteristic $p,$ consider the embedding of the twisted cotangent space $T^*_{X_s'}$ in $\mathbb{P}^{2n}_{k(s)}$ associated to the Rees ring of the  
filtration of $k(s)[x_1^p, \dots, x_n^p, \partial_1^p, \dots, \partial_n^p]$ by the degree of $k(s)$-polynomials in the variables $\{x_1^p, \dots, x_n^p, \partial_1^p, \dots, \partial_n^p\},$ each of which is of degree $1$ (see Definition \ref{defi: filtration} and Proposition \ref{prop: Rees}). Then there is a dense open subset $U \subset S$ such that, for all closed points $s$ of $S$ and for every generic point $z$ of an irreducible component of the $p$-support of $M_s$ of closure $\overline{\{z\}}$ in $\mathbb{P}^{2n}_{k(s)},$ one has $$\deg(\overline{\{z\}})\leq e(M_\mu), \ \mathrm{rk}_z(M_s)\leq e(M_\mu)p^n,$$ where $e(M_\mu)$ is the multiplicity of $M_\mu$ for the Bernstein filtration of the Weyl algebra $A_n(k(\mu))$ (see Definition \ref{defi: multiplicity}), $\mathrm{rk}_z(M_s):=\dim_{k(z)}({(F_*M_s)}_z\otimes k(z))$ and $\deg(\overline{\{z\}})$ is the degree of the reduced closure of the image of $z$ in $\mathbb{P}^{2n}_{k(s)}.$

\item \label{item: splitting} (Theorem \ref{thm:splitting}) There is an open dense subset $U\subset S$ such that for all closed points $s$ of $U,$ the Azumaya algebra $F_*D_{X_s}$ splits on the regular locus of the $p$-support of $M_s,$ where $X_s \xrightarrow{F} X_s'$ is the Frobenius.

\item \label{item: main thm} (Theorem \ref{thm:main thm}) There is an open dense subset $U\subset S$ such that for all closed points $s$ of $U,$ the $p$-support of $M_s$ is a Lagrangian subvariety of $T^*_{X_s'}.$  

\end{enumerate} 

Finally as a corollary of (\ref{item: main thm}), which may be seen as the main result of the paper, we give a new proof that the singular support of a holonomic $\mathcal{D}$-module is a Lagrangian subvariety of the cotangent space, by reduction modulo $p,$ see Corollary \ref{cor: Gabber}.    

The statements (\ref{item: purity}) and (\ref{item: main thm}) are the first essential steps in a comprehensive program of study of the geometry of the $p$-supports of a holonomic $\mathcal{D}$-module as $p$ varies,  outlined by Kontsevich in \cite{Kontsevich}.  

The geometry of the $p$-supports of a given module is very rich indeed. They need neither be conical nor come by reduction modulo $p$ from an invariant defined over $\Spec{\mathbb{Z}},$ and are closely related to the $p$-curvatures of the $\mathcal{D}$-module. Let us illustrate this with a couple examples:

\begin{enumerate}

\item
Let $S=\Spec{\mathbb{Z}},$ $X=\mathbb{A}^n_{S}$ and let $M$ be the finitely generated left $D_{X/S}$-module corresponding to the integrable connection $\nabla=d+dg$ on $\mathcal{O}_X,$ where $g$ is a global section of $\mathcal{O}_X.$ From the identity $(\partial_i+\partial g/\partial x_i)^p=(\partial_i)^p+(\partial g/\partial x_i)^p$ in $A_n(\mathbb{Z}/p\mathbb{Z})$ \cite[5.2.4]{Katz} follows that the $p$-support of $\mathbb{Z}/p\mathbb{Z}\otimes_{\mathbb{Z}}M$ $\subset T^*_{(\mathbb{A}^{n}_{\mathbb{Z}/p\mathbb{Z}})'}=T^*_{\mathbb{A}^n_{\mathbb{Z}/p\mathbb{Z}}}$ is the graph of $dg$ modulo $p.$ 
Thus the $p$-supports are not necessarily conical. 

\item
Let $\mathbb{Z}[\lambda]$ be the subring of $\mathbb{C}$ generated by $\lambda\in\mathbb{C}$ and $S= \Spec\mathbb{Z}[\lambda].$ Let $X=\Spec \mathbb{Z}[\lambda][x,x^{-1}]=\mathbb{A}^1_S-\{0\}\subset \mathbb{A}^1_S=\Spec \mathbb{Z}[\lambda][x].$ Consider the finitely generated left $D_{X/S}$-module $M$ corresponding to the integrable $S$-connection $\nabla=d+\lambda~dx/x$ on $\mathcal{O}_X.$ The identity $(x\partial)^p=x^p\partial^p+x\partial$ in $A_1(\mathbb{Z}/p\mathbb{Z})$ \cite[Lemma 1]{Hochschild} implies that for each closed point $s$ of $S$ of  
positive characteristic $p,$ the $p$-support of $M_s,$  
$p$-supp($M_s$) $\subset T^*_{X'_s}\subset T^*_{\mathbb{A}^{1'}_{k(s)}}=T^*_{\mathbb{A}^1_{k(s)}}$ is given by the equation $xy=\lambda^p-\lambda$ ($\mathrm{mod}$ $p$), where $y$ is the global section of $\mathcal{O}_{T^*_{\mathbb{A}^1_{k(s)}}}$ corresponding to $dx.$ Thus if $\lambda$ is not rational, then the $p$-supports depend nontrivially on $p.$

\end{enumerate}

\subsection{Survey of the proofs}

Let us now comment on the proofs of (\ref{item: purity})-(\ref{item: main thm}) above. This will also serve as a description of the contents of the paper. 

Since pure coherent sheaves have equidimensional supports, we prove (\ref{item: purity}) by showing that $F_*M_s$ is a pure coherent sheaf of dimension $n$ on $T^*_{X'_s},$ for all closed points $s$ of a dense open subset $S_1$ of $S.$ For which we use the following criterion, see Theorem \ref{thm:ext and coherent purity}:  $\forall l\neq n$ and $\forall s\in S_1,$ $$\mathcal{E}\mathnormal{xt}^l_{\mathcal{O}_{T^*_{X'_s}}}(\mathcal{E}\mathnormal{xt}^l_{\mathcal{O}_{T^*_{X'_s}}}(F_*{M_s},\mathcal{O}_{T^*_{X'_s}}),\mathcal{O}_{T^*_{X'_s}})=0.$$ The vanishing of $\mathcal{E}\mathnormal{xt}^l_{D_{X_s}}(\mathcal{E}\mathnormal{xt}^l_{D_{X_s}}(M_s, D_{X_s}), D_{X_s}), \forall l\neq n$ follows from the well-known duality property of holonomic $D_{X_\mu}$-modules by specialization of $M_\mu$ (Proposition \ref{prop:holonomic fibers}). One then checks the criterion using that $F_*D_{X_s}$ is an Azumaya algebra on $T^*_{X'_s}$ (Proposition \ref{prop:pure Azumaya}). This concludes the proof of (\ref{item: purity}). Note that we also prove along the way that if $Y$ is a smooth variety over a perfect field of positive characteristic, then the dimensions of a coherent $D_Y$-module as a $D_Y$-module and as a coherent module over the center of $D_Y$ are equal (Proposition \ref{prop:pure Azumaya}). 

We show in \ref{subsec: reduction} that the proof of the main theorem (\ref{item: main thm}) reduces to $X/S = \mathbb{A}^n_S/S.$ In this case we prove that for all closed points s of a dense open subset $S_2$ of $S_1,$ each irreducible component of the $p$-support $p\text{-}supp(M_s)$ of $M_s$ contains a dense smooth open subset $U$ which is a certain specialization  
of the complement $Y$ of a divisor with normal crossings $D$ in a smooth projective variety over a field of characteristic zero. More precisely there are an irreducible scheme $H$ of generic point $\gamma$ of characteristic zero, a smooth $H$-scheme $\Upsilon$ and a closed point $h\in H$ such that $\Upsilon_h=U$ and $\Upsilon_\gamma=Y.$ 
Moreover we show that there is a differential 1-form $\nu$ on $Y$ such that $d\nu$ has logarithmic poles along $D$ and the restriction of the symplectic form to $U$ is the specialization at $h$  
of $d\nu.$ But by Hodge theory (\cite[Corollaire 3.2.14]{Hodge II}) such a $d\nu$ has to vanish. Finally by our choice of $\Upsilon$ this implies that the symplectic form vanishes on $U.$ Thus the symplectic form vanishes on a dense open subset of the $p$-support and since the $p$-support is equidimensional of dimension $n$ by (\ref{item: purity}), it is a Lagrangian subvariety of $T^*_{(\mathbb{A}_{k(s)}^{n})'}.$ 

Let us now give more details of how this is achieved.
The starting point of the proof are the estimates of (\ref{item: bound}). We verify them as follows.  
Let $\Gamma$ be a good filtration of $M,$ see Definition \ref{defi: good filtration}. Then for all closed points $s$ of a dense open subset $S_3$ of $S_2,$ $\Gamma$ specializes to a good filtration $\Gamma_s$ of $M_s$ such that its Hilbert polynomial is equal to that of the good filtration $\Gamma_\mu$ of the (characteristic zero) $D_{X_\mu}$-module $M_\mu$ (Lemma \ref{lmm:fibers dimension and multiplicity}). In particular, the degree of this Hilbert polynomial is the dimension $m$ of the singular support of $M_\mu$ and its leading coefficient is $m!e(M_\mu).$ We then show that to a good filtration $\Gamma_s$ of $M_s$ is associated a good filtration $p\Gamma_s$ of $M_s$ as a module over the center of $A_n(k(s))$ with the Bernstein filtration (Lemma \ref{lmm:pGamma is good}). The Hilbert polynomial of $p\Gamma_s$ is of degree $m$ and its leading coefficient is $m!e(M_\mu)p^m$ (Proposition \ref{prop:Hilbert polynomial of pGamma}). Thus the coherent sheaf on $\mathbb{P}_{k(s)}^{2n}$ corresponding to the Rees module of $p\Gamma_s$ is of dimension $m$ and degree $e(M_\mu)p^m.$ It then follows from intersection theory that $\Sigma_z \mathrm{rk}_z(M_s)\deg(\overline{\{z\}}) \leq e(M_\mu)p^m,$ where the sum is over the generic points of the (top-dimensional) irreducible components of the $p$-support of $M_s$ (Proposition \ref{prop:mu, ranks and degrees}). This gives the estimate of the ranks since $m = n$ by assumption that $M_\mu$ is holonomic. But $F_*M_s$ is a module over an Azumaya algebra of rank $p^{2n},$ thus by Morita theory the ranks $\mathrm{rk}_z(M_s)$ are divisible by $p^n.$ This gives the estimate of the degrees and concludes the proof of (\ref{item: bound}). Note that it also proves an estimate of the number of irreducible components of $p\text{-}supp(M_s),$ namely $\leq e(M_\mu).$ Let us remark that if we relaxed the holonomicity assumption on $M_\mu$ we would obtain an estimate $\leq e(M_\mu)p^{m-n}$ of the degrees of the top-dimensional irreducible components and it is only in the holonomic case $m=n$ that we get a bound independent of the prime $p.$ This is crucial in what follows.

We  
consider an open embedding $T_{\mathbb{A}_{S_3}^n}^* \hookrightarrow \mathbb{P}_{S_3}^{2n}$ which specializes at every closed point $s$ of $S_3$ to the open embedding of $T^*_{({\mathbb{A}_{k(s)}^{n}})'}= T^*_{\mathbb{A}_{k(s)}^n}$ in $\mathbb{P}^{2n}_{k(s)}$ from (\ref{item: bound}). Let $\mathcal{H}$ be the Hilbert scheme of closed subschemes of $\mathbb{P}_{S_3}^{2n}$ of dimension $n$ and degree $\leq e(M_\mu).$ Provided the estimate of the degrees of the irreducible components of the $p$-support of (\ref{item: bound}), all of dimension $n$ by (\ref{item: purity}), we see that for all closed points $s\in S_3,$ the closure in $\mathbb{P}^{2n}_{k(s)}$ of each irreducible component of $p\text{-}supp(M_s)$ corresponds to a closed point of $\mathcal{H}.$
This allows us to use $\mathcal{H}$ to show that there is an integer $N_1>0$ such that for $S_4= S_3[\frac{1}{N_1}]$ the following is satisfied. 
For each irreducible component $Z$ of $p\text{-}supp(M_s)$ for $s$ a closed point of $S_4,$ there is an open embedding $U\xhookrightarrow{j} \overline{U}$ where $U\subset Z$ is a dense smooth open subset and $\overline{U}$ is a smooth projective variety over $k(s)$ such that the complement of $j$ is a divisor with normal crossings $D.$ Indeed let $\mathcal{Z}_{\mathcal{H}} \subset \mathbb{P}_{\mathcal{H}}^{2n}$ be the universal closed subscheme. By the resolution of singularities in characteristic zero, neglecting finitely many positive characteristics $\{p_1,\dots, p_r\}$ with $p_1\dots p_r = N_1,$ we have a finite partition of the open subset of the Hilbert scheme $\mathcal{H}[\frac{1}{N_1}]$ into irreducible subschemes generically of characteristic zero $(\mathcal{H}_i)$  
and above each $\mathcal{H}_i$ an open embedding $\mathcal{Y}\hookrightarrow \overline{\mathcal{Y}}$ with complement $\Delta$ a divisor with normal crossings relative to $\mathcal{H}_i,$ where $\mathcal{Y}$ is the smooth locus of the intersection $\mathcal{Z}_{\mathcal{H}}|_{\mathcal{H}_i}\cap T_{\mathbb{A}_{\mathcal{H}_i}^n}^*\subset \mathbb{P}_{\mathcal{H}_i}^{2n}$ and $\overline{\mathcal{Y}}$ is a smooth projective $\mathcal{H}_i$-scheme (Proposition \ref{prop:compactification}). Let $z$ be the closed point of the Hilbert scheme $\mathcal{H}$ corresponding to the closure in $\mathbb{P}_{k(s)}^{2n}$ of the irreducible component $Z$ of $p\text{-}supp(M_s)$ and suppose that $z\in \mathcal{H}_i.$ Then we set $(H, h, \gamma, \Upsilon, Y, U\xhookrightarrow{j} \overline{U}, D, \nu)= (\mathcal{H}_i, z, \mu_i, \mathcal{Y}, \mathcal{Y}_{\mu_i}, \mathcal{Y}_z\hookrightarrow \overline{\mathcal{Y}}_z, \Delta_{\mu_i}, \theta_{\mathcal{H}}|_{\mathcal{Y}_{\mu_i}})$ where $\mu_i$ is the generic point of $\mathcal{H}_i$ 
and $\theta_{\mathcal{H}}$ is the canonical 1-form on $T_{\mathbb{A}_{\mathcal{H}}^n}^*.$

We have that the restriction $\omega|_U$ to $U$ of the symplectic form $\omega$ on $T^*_{({\mathbb{A}_{k(s)}^{n}})'}$ is equal to the restriction of the exterior derivative $d\theta_{\mathcal{H}}$ to the fiber $\mathcal{Y}_z.$ Moreover the partition of $\mathcal{H}$ into $(\mathcal{H}_i)$ is such that if $(d\theta_{\mathcal{H}}|_{\mathcal{Y}})_{\mu_i}=0,$ then $d\theta_{\mathcal{H}}|_{\mathcal{Y}}=0$ and if there is a closed point $z$ of $\mathcal{H}_i$ such that  $(d\theta_{\mathcal{H}}|_{\mathcal{Y}})_z$ has logarithmic poles along $\Delta_z,$ then $(d\theta_{\mathcal{H}}|_{\mathcal{Y}})_{\mu_i}$ has logarithmic poles along $\Delta_{\mu_i}$ (Proposition \ref{prop:compactification}). Thus to prove the main theorem (\ref{item: main thm}) it is enough to show that there exists a dense open subset $S_5\subset S_4$ such that for all closed points $s\in S_5$ and for each irreducible component $Z$ of $p\text{-}supp(M_s),$ the restriction of the symplectic form $\omega|_U$ has logarithmic poles along $D,$ where $U$ and $D$ are constructed above. To do so we first show that $\theta_{\mathcal{H}}|_U$ is locally in the image of the $p$-curvature operator $W^*-C_U,$ see Definition \ref{defi: p-curvature}. And then that if a 1-form $\eta$ is locally in the image of the $p$-curvature operator (and has poles along $D$ of order at most $p-1,$ which $\theta_{\mathcal{H}}|_U$ has for all large enough characteristic $p$), then it has logarithmic poles along $D$ (Proposition \ref{prop:log poles}).  

The proof that $\theta_{\mathcal{H}}|_U$ is locally in the image of the $p$-curvature operator goes as follows. There is a morphism $\Omega_U^1 \xrightarrow{\phi_U} Br(U)$ with values in the Brauer group (arising from the $p$-curvature exact sequence, see Definition \ref{defi: phi}) such that $\phi_U(\theta_{\mathcal{H}}|_U)$ is the class of the Azumaya algebra $F_*D_{\mathbb{A}_{k(s)}^n}|_U$ (Proposition \ref{prop:canonical form maps to algebra of differential operators}) and the kernel of $\phi_U$ is the space of 1-forms locally in the image of the $p$-curvature operator (Proposition \ref{prop:kernel of phi}). We thus are left to show that $F_*D_{\mathbb{A}_{k(s)}^n}$ splits on the regular locus of each irreducible component of $p\text{-}supp(M_s).$ This amounts to the splitting of the central simple algebra $\mathcal{A}_z:= (F_*D_{\mathbb{A}_{k(s)}^n})_z\otimes k(z)$ for each irreducible component $Z$ of $p\text{-}supp(M_s)$ of generic point $z.$ Note that this is equivalent to (\ref{item: splitting}) in the case $X/S= \mathbb{A}_S^n/S.$ (From which the general case follows, see the proof of Theorem \ref{thm:splitting}.) But by the estimate of the ranks in (\ref{item: bound}) we have that this central simple algebra has a representation of dimension $\leq e(M_\mu)p^n$ and thus that $e[\mathcal{A}_z]=0$ for some $e\leq e(M_\mu)$ (Lemma \ref{lmm:Azumaya acting on a vector bundle}). Moreover $\mathcal{A}_z$ is of rank $p^{2n},$ thus $p^n[\mathcal{A}_z]=0.$ Hence for $p$ large enough, since $e$ and $p^n$ are coprime, we have $[\mathcal{A}_z]=0.$ We thus see that there exists an integer $N_2>0$ such that $S_5=S_4[\frac{1}{N_2}]$ has the required properties. This concludes the proofs of (\ref{item: splitting}) and (\ref{item: main thm}). 

Let us mention that provided (\ref{item: purity})-(\ref{item: splitting}), one may also prove (\ref{item: main thm}) by adapting the arguments of \cite{Gabber}, as explained in \cite{vdB}. Our original approach presented here has, among other things, the advantage of being independent of \cite{Gabber}, providing in particular a new, more geometric insight into the classical involutivity theorem. 

Note finally that for the reader's convenience we have included an appendix on the (algebraic) symplectic geometry of the cotangent space. 

\subsection{Acknowledgments} 

This paper is a streamlined version of my 2010 Orsay PhD thesis, made under the direction of Maxim Kontsevich. This work owes him a lot. In fact, the systematic use of the Azumaya property excepted, all the main characters involved here were present in \cite{Kontsevich}. He also suggested that ``the restriction of the canonical form to the $p$-support should have logarithmic poles". It is my pleasure to thank him here. I have also had wonderful conversations on questions related to this work with Joseph Bernstein, for which I am very grateful. The influence of Ofer Gabber on the general topics discussed here is obvious. He also provided me with very useful references. Many thanks go to him. 

I would like to thank Pierre Berthelot for his comments on a previous version of this manuscript. The feedback and support from Konstantin Ardakov, Sam Raskin, Simon Riche and Olivier Schiffmann has been quite useful at various stages of the completion of this project. A special thanks goes to them. Finally I am grateful to the referees for their many useful remarks.
The author was partially supported by EPSRC grant EP/L005190/1.

\subsection{Conventions} \label{subsec:general assumptions}

Schemes are assumed to be Noetherian, positive characteristics to be nonzero and morphisms of algebras to preserve the identity element. For $X/S$ a scheme, $s\in S$ and $\Spec k(s) \xrightarrow{i} S$ the corresponding point, we let $X_s/\Spec k(s)$ be the fiber of $X/S$ at $s,$ that is the base-change of $X/S$ by $i.$ If $M$ is a coherent left $D_{X/S}$-module, we denote by $M_s$ the left $D_{X_s/\Spec k(s)}$-module $k(s)\otimes_{\mathcal{O}_S}M,$ the restriction of $M$ to the fiber $X_s/\Spec k(s).$ When there is no risk of confusion, we denote the fiber by $X_s$ (instead of $X_s/\Spec k(s)$) and $D_{X_s/\Spec k(s)}$ by $D_{X_s}.$ As a general rule, if there is no risk of confusion we omit the base scheme $S$ from the notation if $S$ is the spectrum of a field. Local coordinates of a smooth scheme $X/S$ mean local \'{e}tale relative coordinates in the neighborhood of a closed point of $X.$ If $Y$ is a scheme over a field $k$ of positive characteristic, we denote $Y'$ its base-change by the Frobenius endomorphism of $k.$ A reduced scheme over a field is called a variety. As a rule we define notions and state results for left modules, we often omit to mention that they easily adapt to right modules. 

\section{Statement of the main result and general reductions in its proof} \label{sec:statement}

\subsection{Preliminary definitions and notations}
\label{subsec:p-support}
 
Let $S$ be a scheme, $X$ be a smooth $S$-scheme of relative dimension $n$ and let $T_{X/S}$ be the tangent sheaf. 
\begin{defi} The sheaf of crystalline differential operators $D_{X/S}$ on $X/S$ is the enveloping algebra $U_{\mathcal{O}_X}(T_{X/S})$ of the Lie algebroid $(T_{X/S}, [-,-]),$ where $[-,-]$ is the Lie bracket on $T_{X/S}.$\end{defi}

Thus $D_{X/S}$ is generated by the structure sheaf $\mathcal{O}_X$ and the tangent sheaf $T_{X/S}$, subject to relations $f.\partial=f\partial,\ \partial.f-f.\partial=\partial(f)$ and $\partial.\partial'-\partial'.\partial=[\partial,\partial'],$ for all $f$ (resp. $\partial, \partial'$) local sections of $\mathcal{O}_X$ (resp. $T_{X/S}$).  
Note that the formation of $D_{X/S}$ commutes with base-change $S'\to S.$ 
Moreover if $S$ is the spectrum of a field of characteristic $0,$ then $D_{X/S}$ is the usual sheaf of algebraic differential operators on $X.$ When $S$ is the spectrum of a field, we often omit the base $S$ from the notations.

We now briefly discuss the coherence of $D_{X/S}.$ Left multiplication by $\mathcal{O}_X$ makes $D_{X/S}$ into an $\mathcal{O}_X$-module and one sees in local coordinates that $D_{X/S}$ is quasi-coherent. Further using local coordinates, one easily checks the following: 

\begin{prop} \label{prop:filtration D}

The sheaf of rings $D_{X/S}$ has a natural  
filtration $D_{X/S}= \bigcup_{m\geq 0} D_{X/S,\leq m}$, defined by $D_{X/S,\leq 0}:=\mathcal{O}_X$ and $D_{X/S,\leq m+1}:=T_{X/S}.D_{X/S,\leq m}+D_{X/S,\leq m}$, whose associated graded sheaf of rings $grD_{X/S}$ is canonically isomorphic to $\mathcal{O}_{T^*_{X/S}}$, the structure sheaf of the cotangent space of $X/S$.

\end{prop}

Therefore by \cite[Corollaires 2.2.5 and 3.1.2]{BerthI}, $D_{X/S}$ is a sheaf of coherent Noetherian rings. By \cite[Proposition 3.1.3]{BerthI}, coherent $D_{X/S}$-modules have the following properties:

\begin{prop} \label{prop:coherentD}
\begin{enumerate}
\item
A left $D_{X/S}$-module is coherent  
if and only if it is quasi-coherent as an $\mathcal{O}_X$-module and its module of sections over any open of an affine covering is a finitely generated left module over the ring of sections of $D_{X/S}.$  
\item
Assume that X is affine. Then the functor of global sections is an equivalence from the category of coherent left $D_{X/S}$-modules to the category of finitely generated left modules over the global sections of $D_{X/S}.$  
\end{enumerate}
\end{prop}

Let $Y$ be a smooth variety of pure dimension $n$ over a perfect field $k$ of positive characteristic $p.$ We denote the relative Frobenius morphism by $Y \xrightarrow{F} Y'.$ Let $D_Y:= D_{Y/\Spec(k)}$ and $T_Y:= T_{Y/\Spec(k)}.$ Recall that by \cite[Lemma 1.3.2]{BMR}, the $\mathcal{O}_{Y'}$-linear map $T_{Y'} \xrightarrow{c'} F_*D_Y$ sending $\partial$ to $\partial^p-\partial^{[p]},$ where $\partial$ is a local section of $T_{Y'}$ and $\partial^{[p]}$ is its $p$-th power in $T_{Y'},$ lands in the center and induces an isomorphism $$\mathcal{O}_{T^*_{Y'}} \xrightarrow{c} F_*Z(D_Y),$$ with $Z(D_Y)$ the center of $D_Y.$ We will thus consider $F_*D_Y$ as an $\mathcal{O}_{T^*_{Y'}}$-algebra. Furthermore, by \cite[Theorem 2.2.3]{BMR}, $F_*D_Y$ is an Azumaya algebra over $T^*_{Y'}.$ In particular, it is a coherent sheaf on $T^*_{Y'}.$ Hence if $M$ be a coherent left $D_Y$-module, then $F_*M$ is a coherent $\mathcal{O}_{T^*_{Y'}}$-module.

We now introduce our main object of study:

\begin{defi} \label{defi: p-supp}
Let $M$ be a coherent left $D_Y$-module. The $p$-support of $M$ is the support of the coherent $\mathcal{O}_{T^*_{Y'}}$-module $\mathcal{M}$ deduced from the direct image $F_*M,$ using the isomorphism $c.$ It is a closed subset $p$-supp($M$) of $T^*_{Y'}$, which we endow with its reduced subscheme structure. 

\end{defi}

\begin{rem} \label{rem:localization p-support}

Note  
that the $p$-support commutes with \'{e}tale localization on $Y.$
\end{rem}

\begin{rem}  
The schematic support of $\mathcal{M}$ is a not necessarily reduced subscheme of $T^*_{Y'},$ refining the $p$-support of $M.$ We do not explore this notion further here. 

\end{rem}

\subsection{The statement} \label{subsec:The statement}
If $X$ is an $S$-scheme, $M$ a left $D_{X/S}$-module and $s$ a point of $S,$ we denote by $X_s$ the fiber of $X$ at $s$ and by $M_s$ the left $D_{X_s}$-module deduced from $M$ by base-change. We refer to the appendix for the definitions of symplectic form $\omega$ on the cotangent space (Definition \ref{defi: symplectic form}) and Lagrangian subvariety (Definition \ref{defi: Lagrangian}). Our main result is the following:

\begin{thm} \label{thm:main thm}

Let $S$ be an integral scheme dominant and of finite type over $\mathbb{Z},$ of generic point $\mu.$ Let $X$ be a smooth $S$-scheme of relative dimension $n$ and let $M$ be a coherent left $D_{X/S}$-module. 
Suppose that $M_\mu$ is a nonzero holonomic left $D_{X_\mu}$-module, then there is a dense open subset $U$ of $S$ such that the $p$-support of $M_u$ is a Lagrangian subvariety of $(T^*_{X'_u},\omega_{X'_u}),$ for all closed points $u$ of $U.$

\end{thm}

The proof occupies most of the paper and is concluded in Subsection \ref{subsec:conclusion}.

\subsection{First reductions}

Here we carry out some standard reductions.
It is also convenient to consider the case of  
the fiber of $M$ at the generic point $\mu$ of $S$ being zero. We put these into two remarks:

\begin{rem} \label{rem:reduction to affine}
	The conclusion of Theorem \ref{thm:main thm} depends on $S$ only up to restricting to a dense open subset, and so do its hypotheses. Moreover, the assertion is Zariski-local (even \'{e}tale-local) on $X.$ Indeed Lagrangianity is local and so is the $p$-support, as in Remark \ref{rem:localization p-support}. And the hypotheses are stable by restriction to open coverings. Hence to prove the main theorem we may further assume that $S$ is affine, regular  
	and that $X$ is regular,  
	affine and integral. 
	
\end{rem}

\begin{rem} \label{rem:generic fiber zero}
	
	If the fiber of $M$ at the generic point $\mu$ of $S$ is zero,  
	then there is a dense open subset $U$ of $S$ such that $M|_{U}=0.$
	
\end{rem}

\begin{proof}
 
Indeed one may assume that $X$ and $S$ are affine and thus consider a left  
module over the ring of global sections of $D_{X/S}.$ By the hypotheses, this module has a finite generating family $\{m_1,...,m_l\}$ and each $m_i$ is annihilated by a nonzero global section $r_i$ of $\mathcal{O}_S.$ Since $\mathcal{O}_S$ acts through the center of $D_{X/S},$ the open subset of $S$ determined by the product of these global sections fulfills the statement.
\end{proof}

\subsection{Reduction to $\mathbb{A}^n$} \label{subsec: reduction}
Here we show that the proof of the main theorem reduces to $X/S=\mathbb{A}^n_S.$ To do so we use the direct image of $D_{X/S}$-modules and the general result on the dimension of $p$-supports (Theorem \ref{thm:purity}), proved independently below.

\begin{prop} \label{prop:reduction to affine space} To prove Theorem \ref{thm:main thm}, it is sufficient to suppose that $X/S=\mathbb{A}^n_S.$  
\end{prop}
\begin{proof}
By Remark \ref{rem:reduction to affine}, one may suppose that $X/S$ is smooth of relative dimension $n$ and that $X$ and $S$ are affine. Hence there is a closed immersion $X \overset{f}\hookrightarrow \mathbb{A}^m_S$ over $S,$ for some $m\geq0.$ 

Let $M$ be a left $D_{X/S}$-module as in the statement of Theorem \ref{thm:main thm} and let $f_+M$ be its direct image, see \cite[2.4.1]{BerthIntro} for the definition. It is easy to see that one has the classical description of the transfer bimodule, as in \cite[(4) p.259]{Borel}. Hence the latter is flat over $D_{X/S}.$ In particular, the direct image $f_+M$ is supported in degree $0$ and $H^0f_+M$ is a coherent left $D_{\mathbb{A}^m_S/S}$-module.  
Furthermore, it follows directly from the definition that the formation of the transfer bimodule commutes with base-change $S'\to S.$  
Hence $H^0f_+$ commutes with base-change $S'\to S.$  
Thus the generic fiber $(H^0f_+M)_\mu= H^0{f_\mu}_+M_\mu$ is nonzero and is holonomic by preservation of holonomicity under direct images, where $\mu$ is the generic point of $S$ and $X_\mu \overset{f_\mu}\hookrightarrow \mathbb{A}^m_{k(\mu)}$ is the induced closed immersion. Finally, for all closed points $s\in S,$ it is also an immediate consequence of the description of the transfer bimodule as in \cite[(4) p.259]{Borel} that $$p\textit{-supp}(H^0{f_s}_+M_s)= ({f_s'})_{\pi}\circ({f_s'})_d^{-1}(p\text{-}supp(M_s)),$$ where $f_s'$ is the base-change of $f_s$ by the Frobenius and we have used the notation of the appendix for the maps in the cotangent diagram of $f_s'.$

Thus by Lemma \ref{lmm:isotropy}, to prove that the symplectic form vanishes on the regular locus of $p\text{-}supp(M_s)$ it is enough to prove the corresponding symplectic form vanishes on the regular locus of $p\textit{-supp}(H^0{f_s}_+M_s).$ This concludes the proof of the proposition since the part of the theorem concerning dimensions is Theorem \ref{thm:purity}, proved independently below.\end{proof}

\section{Dimension of the $p$-supports} \label{sec:dimension of the p-supports}

As outlined in the introduction,  
the main theorem splits into an assertion about the dimension of the $p$-support and one about the vanishing of the symplectic form on the regular locus of the $p$-support. We start by considering the former.

\subsection{Statement}

\begin{thm} \label{thm:purity}
Let $S$ be an integral scheme dominant and of finite type over $\mathbb{Z},$ of generic point $\mu.$ Let $X$ be a smooth $S$-scheme of relative dimension $n$ and let $M$ be a coherent left $D_{X/S}$-module. 
Suppose that $M_\mu$ is a nonzero holonomic left $D_{X_\mu}$-module, then there is a dense open subset $U$ of $S$ such that the $p$-support of the fiber of $M$ at each closed point $u$ of $U$ is equidimensional of dimension $n = \dim X.$

\end{thm}

The proof is contained in the subsection \ref{subsec:equidim of the p-supports}. 
In view of Remark \ref{rem:reduction to affine}, we may and shall assume that $S$ and $X$ are regular, integral and affine. 

\subsection{Pure coherent sheaves} \label{subsec:pure sheaves}

Recall that the \textit{(co)dimension} of a coherent sheaf is the (co)dimension of its support and let us call a coherent sheaf \textit{equidimensional} if its support is equidimensional. There is a strengthening of equidimensionality which has a very convenient interpretation in terms of duality theory. Namely, let $Y$ be an affine scheme.

\begin{defi}  A coherent sheaf on $Y$ is \textit{pure} if all its nonzero coherent subsheaves are of the same dimension. 

\end{defi}

It is easily seen to imply equidimensionality:

\begin{prop} \label{prop:pure equidimensional}

A coherent sheaf $\mathcal{F}$ on $Y$ is pure if and only if all its associated points $y \in Ass(\mathcal{F})$ are of the same dimension. In particular, a pure coherent sheaf on $Y$ is equidimensional.

\end{prop}

\begin{proof} By \cite[Proposition 3.1.2]{EGAIV} a prime ideal $\mathfrak{p}$ corresponding to an associated point of $\mathcal{F}$ is associated to the module of global sections of $\mathcal{F}.$ Thus $\mathfrak{p}$ is the annihilator of a global section of $\mathcal{F}.$ The only if part immediately follows.

Suppose that the dimension of every associated point of $\mathcal{F}$ is $d$ and let $\mathcal{F}'\subset\mathcal{F}$ be a nonzero coherent subsheaf. Then the dimension of $\mathcal{F}'$ is $d.$ Indeed the associated points of $\mathcal{F}'$ contain the generic points of its support, by \cite[Corollaire 3.1.4]{EGAIV}, and are associated to $\mathcal{F}$ by \cite[Proposition 3.1.7 (i)]{EGAIV}. This concludes the proof of the proposition.
\end{proof}

Here is the interpretation in terms of duality theory.

\begin{thm} \label{thm:ext and coherent purity}

Suppose that $Y$ is regular and equidimensional. A coherent sheaf $\mathcal{F}$ on $Y$ is pure if and only if there is a nonnegative integer $c$ such that $$\mathcal{E}\mathnormal{xt}_{\mathcal{O}_Y}^l(\mathcal{E}\mathnormal{xt}_{\mathcal{O}_Y}^l(\mathcal{F},\mathcal{O}_Y),\mathcal{O}_Y)=0$$ for all $l \neq c$. If $\mathcal{F}$ is nonzero, then $c$ is its codimension.

\end{thm}

\begin{proof} This is well-known. We refer to the literature. Our definition of purity is equivalent to \cite[A:IV 2.5.]{Bjork2} by \cite[V, 2.2.3]{Borel}. The theorem is then \cite[A:IV 2.6]{Bjork2}, since a regular local ring is Auslander regular by \cite[A:IV 3.4]{Bjork2}. \end{proof}

\subsection{Equidimensionality of the $p$-supports} \label{subsec:equidim of the p-supports}
Here we prove Theorem \ref{thm:purity}. In particular, we use Remark \ref{rem:reduction to affine} and suppose that $X/S$ is smooth of relative dimension $n$ with $S$ and $X$ regular, affine and integral. We start by recalling the notion of good filtration on a coherent $D_{X/S}$-module. Recall the filtration on $D_{X/S}$ from Proposition \ref{prop:filtration D}.
\begin{defi} \label{defi: good filtration}
A good filtration on a coherent left $D_{X/S}$-module is a filtration by coherent $\mathcal{O}_X$-submodules, compatible with the filtration on $D_{X/S},$ which is bounded below and such that the associated graded module over $grD_{X/S} \cong \mathcal{O}_{T^*_{X/S}}$ is coherent. 
\end{defi}
Note that coherent left $D_{X/S}$-modules admit good filtrations by \cite[5.2.3 (iv)]{BerthIntro}. 

We next give a lemma guaranteeing the freeness on $S$ of a module whose associated graded is free. The proof of the last part of the lemma was kindly provided by Michel Van den Bergh.   
 
 \begin{lmm} \label{lmm:graded flatness}

Let $M$ be a left module over a ring $R$ and let $\{M_i\}_{i\in \mathbb{Z}}$ be an exhaustive increasing filtration of $M$ by left $R$-submodules. 
Suppose that there is $i_0 \in \mathbb{Z}$ such that $M_{i_0}=0$ and, $\forall i > i_0,$ the left R-module $M_i/M_{i-1}$ is flat, then $M$ is flat. 
Suppose further that $\forall i, M_i/M_{i-1}$ is free, then $M$ is free. 
\end{lmm}
\begin{proof} By hypothesis, $M_{i_0+1} \cong M_{i_0+1}/(M_{i_0}=0)$ is flat. Moreover, $\forall i \geq i_0+1, M_i/M_{i-1}$ is flat. So, since extensions of flat modules are flat (\cite[\S 2 n$^{\rm o}$5 Proposition 5]{BourbakiACI}), $\forall i \geq i_0+1, M_i$ is flat. Thus $M$ is a union of flat submodules. Hence it is flat by \cite[\S 2 n$^{\rm o}$3 Proposition 2(ii)]{BourbakiACI}. This proves the first assertion. 

Suppose that the $M_i/M_{i-1}$ are free. Then the union, over all $i \geq i_0+1,$ of an arbitrary lift to $M_i$ of a basis of $M_i/M_{i-1}$ is a basis of $M.$ Thus $M$ is free. This finishes the proof of the lemma.
\end{proof}
  
  The following lemma is standard for $\mathcal{O}$-modules. Considering the associated graded allows us to deduce a version for $D_{X/S}$-modules. 
    
\begin{lmm} \label{lmm:ext fibers}

Let $M$ be a coherent left $D_{X/S}$-module. Then there is a dense open subset $U$ of $S$ such that $\forall l\geq0$ and $\forall s \in U$, the canonical map $$(\mathcal{E}\mathnormal{xt}^l_{D_{X/S}}(M,D_{X/S}))_s \to \mathcal{E}\mathnormal{xt}^l_{D_{X_s}}(M_s,D_{X_s})$$ is an isomorphism, where the subscript $s$ denotes the restriction to the fiber.

\end{lmm}
\begin{proof} First of all, there are only finitely many degrees $l$ to consider. Namely, by \cite[A:IV 4.5]{Bjork2}, both target and domain of the above morphism are zero for $l > \dim T^*_{X/S} \geq \dim T^*_{X_s}.$ Indeed, $T^*_{X/S}$ and $T^*_{X_s}$ are the respective spectra of the rings $grD_{X/S}$ and $grD_{X_s},$ which are both regular. Thus guaranteeing the vanishing of $\mathcal{E}\mathnormal{xt}^l$ for $l > \dim T^*_{X/S}.$ Hence it suffices to prove that, $\forall l,$ there is an open $U$ as in the statement.

We thus want to prove a $D_{X/S}$-module version of \cite[Corollaire 9.4.3]{EGAIV}. But the proof of the latter adapts to $D_{X/S}$-modules by considering the associated graded to good filtrations. Indeed, coherent left $D_{X/S}$-modules form an abelian category and the proof of \cite[Proposition 9.4.2]{EGAIV} carries through, using \cite[Proposition A.17]{Kash} and Lemma \ref{lmm:graded flatness} to conclude. \end{proof}

We now obtain that the following well-known consequence of holonomicity spreads from the generic fiber to a neighborhood.

\begin{prop} \label{prop:holonomic fibers}

Let $M$ be a coherent left $D_{X/S}$-module. Suppose that $M_\mu$ is a holonomic left $D_{X_\mu}$-module, for $\mu$ the generic point of $S.$ Then there is a dense open subset $U$ of $S$ such that for all $l \neq n$ and all $s \in U$, $$\mathcal{E}\mathnormal{xt}^l_{D_{X_s}}(M_s,D_{X_s})=0.$$

\end{prop}
\begin{proof} By Lemma \ref{lmm:ext fibers} and  \cite[VI 1.12]{Borel}, $\forall l \neq n,$ the fiber of $\mathcal{E}\mathnormal{xt}^l_{D_{X/S}}(M,D_{X/S})$ at the generic point of $S$ vanishes. Hence by Remark \ref{rem:generic fiber zero} and Lemma \ref{lmm:ext fibers}, $\forall l \neq n,$ there is a dense open subset $U_l$ of $S$ such that for all $s \in U_l,$ $\mathcal{E}\mathnormal{xt}^l_{D_{X_s}}(M_s,D_{X_s})=0.$ Since by the proof of \ref{lmm:ext fibers} there are only finitely many such degrees $l$ to consider, $U:=\bigcap_lU_l$ fulfills the proposition.\end{proof}

Then we use the Azumaya property of the ring of differential operators in positive characteristic (\cite[Theorem 2.2.3]{BMR}) to transfer purity from a $D_{X/S}$-module to its associated coherent sheaf on the twisted cotangent space.  

\begin{prop} \label{prop:pure Azumaya}

Let $Y$ be a smooth equidimensional scheme over a field $k$ of positive characteristic $p,$ let $Y \xrightarrow{F} Y'$ be the relative Frobenius and let $M$ be a coherent left $D_Y$-module. Then, $\forall l\geq0,$ $$\mathcal{E}\mathnormal{xt}^l_{D_Y}(M,D_Y)=0 \text{ if and only if } \mathcal{E}\mathnormal{xt}^l_{\mathcal{O}_{T^*_{Y'}}}(\mathcal{M},\mathcal{O}_{T^*_{Y'}})=0,$$ where $\mathcal{M}:=F_*M$ is endowed with an action of $\mathcal{O}_{T^*_{Y'}}$ as in Definition \ref{defi: p-supp}. 

\end{prop}
\begin{proof} Since $F$ is affine, $\mathcal{E}\mathnormal{xt}^l_{D_Y}(M,D_Y)=0$ if and only if $$0={F}_*\mathcal{E}\mathnormal{xt}^l_{D_Y}(M,D_Y) \cong \mathcal{E}\mathnormal{xt}^l_{{F}_*D_Y}({F}_*M,{F}_*D_Y).$$ 
Set $\mathcal{D}_Y:={F}_*D_Y,$ we thus have $\mathcal{E}\mathnormal{xt}^l_{{F}_*D_Y}({F}_*M,{F}_*D_Y)= \mathcal{E}\mathnormal{xt}^l_{\mathcal{D}_Y}(\mathcal{M},\mathcal{D}_Y).$
Let us show that $$\mathcal{E}\mathnormal{xt}^l_{\mathcal{D}_Y}(\mathcal{M},\mathcal{D}_Y)=0 \text{ if and only if } \mathcal{E}\mathnormal{xt}^l_{\mathcal{O}_{T^*_{Y'}}}(\mathcal{M},\mathcal{O}_{T^*_{Y'}})=0.$$ 

Indeed, both $\mathcal{E}\mathnormal{xt}^l_{\mathcal{D}_Y}(\mathcal{M},\mathcal{D}_Y)$ and $\mathcal{E}\mathnormal{xt}^l_{\mathcal{O}_{T^*_{Y'}}}(\mathcal{M},\mathcal{O}_{T^*_{Y'}}),$ are quasi-coherent sheaves on $T^*_{Y'}.$ Hence their respective vanishings may be checked on a flat covering $\mathcal{U} \xrightarrow{\pi} T^*_{Y'}$ of $T^*_{Y'}$. Since $\mathcal{D}_Y$ is an Azumaya algebra over $\mathcal{O}_{T^*_{Y'}}$ by \cite[Theorem 2.2.3]{BMR}, this covering may be chosen to split $\mathcal{D}_Y.$ That is $(\mathcal{D}_Y)_\mathcal{U}:=\pi^*\mathcal{D}_Y \simeq M_r(\mathcal{O}_{\mathcal{U}})$, the sheaf of $r\times r$ matrices with coefficients in $\mathcal{O}_{\mathcal{U}}.$ 

As is well-known in Morita theory, tensoring with the $(M_r(\mathcal{O}_{\mathcal{U}}),\mathcal{O}_{\mathcal{U}})$-bimodule $\mathcal{O}^r_{\mathcal{U}}$ induces an equivalence between the category of coherent $\mathcal{O}_{\mathcal{U}}$-modules and the category of coherent left $M_r(\mathcal{O}_{\mathcal{U}})$-modules.  
Note that the coherent sheaf $(\mathcal{O}^r_{\mathcal{U}})^\vee$ is sent to $\mathcal{O}^r_{\mathcal{U}} \otimes_{\mathcal{O}_{\mathcal{U}}}(\mathcal{O}^r_{\mathcal{U}})^\vee\cong M_r(\mathcal{O}_{\mathcal{U}})$ by this equivalence. Let $\mathcal{F}$ be a coherent sheaf such that $\mathcal{M}_\mathcal{U}:=\pi^*\mathcal{M}\simeq \mathcal{O}^r_{\mathcal{U}} \otimes_{\mathcal{O}_{\mathcal{U}}}\mathcal{F}$ as coherent left $(\mathcal{D}_Y)_\mathcal{U} \simeq M_r(\mathcal{O}_{\mathcal{U}})$-modules.  
Then, by localization and the above Morita equivalence, $$\pi^*\mathcal{E}\mathnormal{xt}^l_{\mathcal{D}_Y}(\mathcal{M},\mathcal{D}_Y)\simeq \mathcal{E}\mathnormal{xt}^l_{(\mathcal{D}_Y)_\mathcal{U}}(\mathcal{M}_\mathcal{U},(\mathcal{D}_Y)_\mathcal{U}) $$ $$\simeq \mathcal{E}\mathnormal{xt}^l_{M_r(\mathcal{O}_{\mathcal{U}})}(\mathcal{O}^r_{\mathcal{U}} \otimes_{\mathcal{O}_{\mathcal{U}}}\mathcal{F},\mathcal{O}^r_{\mathcal{U}} \otimes_{\mathcal{O}_{\mathcal{U}}}(\mathcal{O}^r_{\mathcal{U}})^\vee)\simeq_{\mathcal{O}_{\mathcal{U}}\textit{-mod}} \mathcal{E}\mathnormal{xt}^l_{\mathcal{O}_{\mathcal{U}}}(\mathcal{F},(\mathcal{O}^r_{\mathcal{U}})^\vee)\text{ vanishes}$$ if and only if, by commutation with finite direct sums, $\mathcal{E}\mathnormal{xt}^l_{\mathcal{O}_{\mathcal{U}}}(\mathcal{F},\mathcal{O}_{\mathcal{U}})\text{ vanishes},$ if and only if $$\mathcal{E}\mathnormal{xt}^l_{\mathcal{O}_{\mathcal{U}}}(\mathcal{O}^r_{\mathcal{U}} \otimes_{\mathcal{O}_{\mathcal{U}}}\mathcal{F},\mathcal{O}_{\mathcal{U}})\simeq \mathcal{E}\mathnormal{xt}^l_{\mathcal{O}_{\mathcal{U}}}(\mathcal{M}_\mathcal{U},\mathcal{O}_{\mathcal{U}})
\simeq \pi^*\mathcal{E}\mathnormal{xt}^l_{\mathcal{O}_{T^*_{Y'}}}(\mathcal{M},\mathcal{O}_{T^*_{Y'}}) \text{ vanishes,}$$ using again commutation with finite direct sums and localization. This concludes the proof of the proposition.\end{proof}

We can now prove the theorem.

\begin{proof}(of Theorem \ref{thm:purity}) Note that if the fiber of $M$ at the generic point of $S$ is nonzero then $M$ is nonzero. Therefore, by generic freeness (\cite[Theorem 14.4]{EisenbudCA}) applied to the associated graded to a good filtration on $M$ and Lemma \ref{lmm:graded flatness}, there is a dense open subset $W$ of $S$ on which $M$ is faithfully flat. Hence $\forall s\in W, M_s\neq 0$ and thus $F_*M_s\neq 0.$ Since $M$ is holonomic on the generic fiber of $S,$ there is a dense open subset $U$ of $W$ such that $\forall l\neq n$ and $\forall s\in U, \mathcal{E}\mathnormal{xt}^l_{D_{X_s}}(M_s,D_{X_s})=0,$ by Proposition \ref{prop:holonomic fibers}. Which, by Proposition \ref{prop:pure Azumaya}, is equivalent to, $\forall l\neq n$ and $\forall s\in U, \mathcal{E}\mathnormal{xt}^l_{\mathcal{O}_{T^*_{X'_s}}}(F_*{M_s},\mathcal{O}_{T^*_{X'_s}})=0.$ In particular, $\forall l\neq n$ and $\forall s\in U,$ $$\mathcal{E}\mathnormal{xt}^l_{\mathcal{O}_{T^*_{X'_s}}}(\mathcal{E}\mathnormal{xt}^l_{\mathcal{O}_{T^*_{X'_s}}}(F_*{M_s},\mathcal{O}_{T^*_{X'_s}}),\mathcal{O}_{T^*_{X'_s}})=0.$$ This implies by Theorem \ref{thm:ext and coherent purity} that $\forall s\in U, F_*{M_s}$ is a pure nonzero coherent $\mathcal{O}_{T^*_{X'_s}}$-module of dimension $n.$ Hence it is equidimensional of dimension $n$ by Proposition \ref{prop:pure equidimensional}. This proves the theorem. \end{proof}

\begin{rem}

The purity of the coherent $\mathcal{O}_{T^*_{X'_s}}$-module $F_*M_s$ guarantees furthermore that it has no embedded associated points.

\end{rem}

\section{Degrees and ranks estimates} \label{sec:bound}

We now consider $D_{\mathbb{A}^n_S/S}$-modules. In addition to the filtration by the order of differential operators (Proposition \ref{prop:filtration D}), $D_{\mathbb{A}^n_S/S}$ is endowed with the Bernstein filtration. The latter has the property that each summand of its associated graded ring is a finitely generated module over $\mathcal{O}(S).$ For a $D_{\mathbb{A}^n_S/S}$-module $M$ whose generic fiber is holonomic, we use the Bernstein filtration to estimate the degree (for a suitable projective embedding) of the $p$-support $p\textit{-supp}(M_s)$ as well as the rank of $F_*M_s$ at the generic point of an irreducible component of its support $p\textit{-supp}(M_s),$ for $s$ a closed point in a dense open subset of $S,$  
see Theorem \ref{thm:bound}. 

\subsection{Bernstein filtration} \label{subsec:Bernstein filtration}
Let $S$ be an affine scheme and let $R$ be its ring of global sections. If we fix coordinates $\{x_1, \dots, x_n\}$ on $\mathbb{A}^n_S,$ then  
the ring of global sections of $D_{\mathbb{A}^n_S/S}$ is isomorphic to the $n$-th Weyl algebra $A_n(R)$ over $R,$$$A_n(R):= R[x_1,...,x_n]\langle\partial_1,...,\partial_n\rangle/ \langle[\partial_i,\partial_j], [\partial_i,x_j]- \delta_{i,j}; \forall 1\leq i,j\leq n \rangle.$$  
\begin{defi} The \textit{Bernstein filtration} $\mathcal{B}$ of $A_n(R)$ is the filtration by the total order in $x$ and $\partial.$ Namely  
$\forall l \in \mathbb{Z},$ $\mathcal{B}_lA_n(R):=\bigoplus_{|\alpha|+|\beta|\leq l} R x^{\alpha}\partial^{\beta},$ where $\alpha, \beta \in \mathbb{Z}_{\geq0}^n$ are multi-indices and we have used the standard notation $x^\alpha:= x_1^{\alpha_1}\dots x_n^{\alpha_n}, \partial^\beta:= \partial_1^{\beta_1}\dots\partial_n^{\beta_n}$ and for a multi-index $\alpha \in \mathbb{Z}_{\geq0}^n, |\alpha|:= \alpha_1+\dots+\alpha_n.$ 

\end{defi}  
\begin{rem} Note that the associated graded ring $gr^{\mathcal{B}}A_n(R)$ is the $R$-algebra of polynomials in the variables $\{x_1,...,x_n,$ $y_1,...,y_n\},$ graded by the order of polynomials. Where, $\forall 1\leq i\leq n, x_i$ (resp. $y_i$) is the class of $x_i$ (resp. $\partial_i$) $\in \mathcal{B}_1A_n(R)/\mathcal{B}_0A_n(R).$  
In particular, $\forall l \in \mathbb{Z}, \mathcal{B}_lA_n(R)/\mathcal{B}_{l-1}A_n(R)$ is a finitely generated free $R$-module.
\end{rem}

We will use the notion of good filtration on a $A_n(R)$-module. 
\begin{defi}
A filtration $\Gamma$ of a left $A_n(R)$-module $M$ is an increasing exhaustive filtration of $M,$ indexed by $\mathbb{Z}$ and compatible with $\mathcal{B}.$ It is said to be a good filtration if it is bounded below and the associated graded module $gr^{\Gamma}M$ is finitely generated over the algebra $gr^{\mathcal{B}}A_n(R)$. \end{defi}
It is easy to see that finitely generated left $A_n(R)$-modules have good filtrations, see e.g. \cite[Ch.1 Proposition 2.7]{Bjork1}. 
Note that if $\Gamma$ is a good filtration of $M,$ then $\forall l\in \mathbb{Z}, \Gamma_lM/\Gamma_{l-1}M$ (and hence $\Gamma_lM$) is a finitely generated $R$-module. 

 Suppose that $R$ is a field $K$. Then a left $A_n(K)$-module has well-defined degree and multiplicity. Indeed, let $M$ be a finitely generated left $A_n(K)$-module and let $\Gamma$ be a good filtration on $M$. Then for $l$ large enough, the function $l\mapsto \dim_K\Gamma_lM$ coincides with a polynomial $\mathcal{H}_{M,\Gamma} \in \mathbb{Q}[t]$ (\cite[Ch.1 Corollary 3.3]{Bjork1}). \begin{defi} \label{defi: multiplicity}
 	Let $d$ (resp. $a_d$) be the degree (resp. the leading coefficient) of $\mathcal{H}_{M,\Gamma}.$ Then $d!a_d$ is a nonnegative integer. The nonnegative integers $d(M):=d$ and $e(M):=d!a_d$ are independent of $\Gamma$ and called the dimension and multiplicity of $M$, respectively (\cite[p.8]{Bjork1}).\end{defi}

Now we look at the behavior of these invariants in a family.

\begin{lmm} \label{lmm:fibers dimension and multiplicity}

Suppose that $R$ is a domain and let $M$ be a finitely generated left $A_n(R)$-module. Then there is a dense open subset $U$ of $S:=\Spec(R)$ such that the functions $s\mapsto d(M_s)$ and $s\mapsto e(M_s)$ are constant on $U$.

\end{lmm}

\begin{proof} Let $\Gamma$ be a good filtration on $M$. Then by generic freeness (\cite[Theorem 14.4]{EisenbudCA}), there is a dense open subset $U$ of $S$ such that $\forall l\in\mathbb{Z}, (\Gamma_lM/\Gamma_{l-1}M)|_U$ is free over $\mathcal{O}(U).$ In particular, $\forall l\in \mathbb{Z}, (\Gamma_lM/\Gamma_{l-1}M)|_U$ is a flat $\mathcal{O}(U)$-module. Hence, $\forall l\in \mathbb{Z}$ and $\forall s\in U, (\Gamma_lM/\Gamma_{l-1}M)_s\cong(\Gamma_lM)_s/(\Gamma_{l-1}M)_s$ and $(\Gamma)_s$ is a good filtration on $M_s$. The lemma follows since, $\forall s\in U$ and $\forall l\in \mathbb{Z}, \dim_{k(s)}(\Gamma_lM)_s=\sum_{i=-\infty}^{i=l} \dim_{k(s)}(\Gamma_iM)_s/(\Gamma_{i-1}M)_s$ and $$\dim_{k(s)}(\Gamma_lM)_s/(\Gamma_{l-1}M)_s=\dim_{k(s)}(\Gamma_lM/\Gamma_{l-1}M)_s$$ is the rank of the free module $\mathcal{O}(U)$-module $\Gamma_lM/\Gamma_{l-1}M|_U.$ Hence $\mathcal{H}_{M_s,\Gamma_s}$ is constant on $U.$  
\end{proof}

\subsection{On the filtrations of the center} \label{subsec:induced filtration over the center}

Let $K$ be a field of positive characteristic $p.$ With the notation of Subsection \ref{subsec:Bernstein filtration}, the center $ZA_n(K)$ of $A_n(K)$ is the algebra of polynomials $K[x_1^p,...,x_n^p,\partial_1^p,...,\partial_n^p].$ \begin{defi} \label{defi: filtration}
 Let $\mathcal{C}$ be the
 filtration on $ZA_n(K)$ derived from the grading $| \bullet |$ of polynomials, where $|x_i^p|= |\partial_j^p|= 1.$
\end{defi} We would like to compare $\mathcal{C}$ with the Bernstein filtration and to do so we will use the classical construction of the Rees ring associated to a filtered ring. 

\begin{defi} \label{defi:Rees} \begin{enumerate}\item The Rees ring  of the filtered ring $(ZA_n(K),\mathcal{C})$ is the graded ring $R_n(\mathcal{C}):=\bigoplus\limits_{i\in\mathbb{N}} \mathcal{C}_iZA_n(K).$ 
\item Let $G$ be an increasing $\mathcal{C}$-compatible filtration of a $ZA_n(K)$-module $M.$ 
Then the Rees module associated with $G$ is $\mathcal{R}(M,G):=\bigoplus\limits_{i\in \mathbb{Z}} G_iM.$ It is a naturally an $R_n(\mathcal{C})$-module. \end{enumerate}\end{defi}

We recall elementary properties of the Rees ring in the following lemma.

\begin{lmm} \label{lmm: Rees ring}
\begin{enumerate}
\item \label{item: iso}The graded algebra morphism $$K[t_0,x_1^p,...,x_n^p,\partial_1^p,...,\partial_n^p] \to R_n(\mathcal{C}):=\bigoplus\limits_{i\in \mathbb{N}} \mathcal{C}_iK[x_1^p,...,x_n^p,\partial_1^p,...,\partial_n^p]$$ $$t_0\mapsto 1,  
x_i^p \mapsto x_i^p, 
\partial_j^p \mapsto \partial_j^p$$
with $t_0$ in degree 1, is an isomorphism. 
\item \label{item: fiber} Using the same notation for $t_0$ and its image under the isomorphism of (\ref{item: iso}), the natural map $R_n(\mathcal{C})/t_0R_n(\mathcal{C}) \to gr^{\mathcal{C}}ZA_n(K)$ is an isomorphism of graded algebras. 
\item \label{item: sum}There is a unique morphism $R_n(\mathcal{C})_{t_0}\to ZA_n(K),$ sending $t_0$ to $1$ and extending the inclusions $\mathcal{C}_iK[x_1^p,...,x_n^p,\partial_1^p,...,\partial_n^p] \subset ZA_n(K).$ 
Its restriction to $R_n(\mathcal{C})_{(t_0)},$ the subring of degree $0$ elements of the graded ring $R_n(\mathcal{C})_{t_0},$ 
is an isomorphism $R_n(\mathcal{C})_{(t_0)}\widetilde{\to} ZA_n(K).$  
\end{enumerate}\end{lmm}

We will use good filtrations of $ZA_n(K)$-modules.
\begin{defi}
A filtration $G$ of a $ZA_n(K)$-module $M$ as in the definition \ref{defi:Rees} is said to be good if the associated Rees module $\mathcal{R}(M,G)$ is a finitely generated $R_n(\mathcal{C})$-module. \end{defi} This implies in particular that $G$ is bounded below. Moreover, one easily sees that a filtration $G$ on $M$ is good if and only if $G$ is bounded below and the associated graded module $gr^GM$ is finitely generated over $gr^{\mathcal{C}}ZA_n(K)$ (\cite[A:III 1.29]{Bjork2}).

Let $M$ be a left $A_n(K)$-module and let $r_*M$ be the module $M$ considered as a $ZA_n(K)$-module. We now introduce a $\mathcal{C}$-filtration of $r_*M$ associated with a $\mathcal{B}$-filtration of $M.$

\begin{defi}
Let $\Gamma$ be a filtration of the left $A_n(K)$-module $M.$ The $\mathcal{C}$-filtration $p\Gamma$ of $r_*M$ is given by $(p\Gamma)_lr_*M:=\Gamma_{pl}M,$ for all integers $l.$ 
\end{defi}

We now want to relate properties of $p\Gamma$ to those of $\Gamma.$ 

\begin{lmm} \label{lmm:pGamma is good}
Let $\Gamma$ be a good filtration of the left $A_n(K)$-module $M,$ then $p\Gamma$ is a good filtration of $r_*M.$ 
\end{lmm}

\begin{proof} It is clear from the definition that the filtration $p\Gamma$ is bounded below, since $\Gamma$ is. 

Let us show that the $gr^{\mathcal{C}}ZA_n(K)$-module $gr^{p\Gamma}r_*M$ is finitely generated. In order to do so, we consider the filtration $\Phi$ of the center induced by the Bernstein filtration.  
Thus, $\forall l\in \mathbb{Z}, \Phi_lZA_n(K):=ZA_n(K)\cap \mathcal{B}_lA_n(K).$ In particular, $x_i^p$ and $\partial_j^p$ are of degree $p$ for the filtration $\Phi,$ for all $i,j.$  
Let $\Phi(\Gamma)$ be the $\Phi$-filtration on $r_*M$ defined by $\Phi(\Gamma)_lr_*M:=\Gamma_{pm}M$, where $pm$ is the greatest integer multiple of $p$ such that $pm\leq l$. Note that there is a $K$-module isomorphism $gr^{p\Gamma}r_*M\to gr^{\Phi(\Gamma)}r_*M$ defined by $$(p\Gamma)_lr_*M/(p\Gamma)_{l-1}r_*M= \Gamma_{pl}M/\Gamma_{p(l-1)}M= \Phi(\Gamma)_{pl}r_*M/\Phi(\Gamma)_{pl-1}r_*M.$$ It is $\rho$-linear, where $\rho$ is  
the isomorphism of $K$-algebras $gr^{\mathcal{C}}ZA_n(K)\to gr^\Phi ZA_n(K)$ satisfying
$$x_i^p + \mathcal{C}_2ZA_n(K) \mapsto x_i^p + \Phi_{p+1}ZA_n(K); \partial_j^p + \mathcal{C}_2ZA_n(K) \mapsto \partial_j^p + \Phi_{p+1}ZA_n(K).$$ Hence $gr^{p\Gamma}r_*M$ is finitely generated over $gr^{\mathcal{C}}ZA_n(K)$ if and only if $gr^{\Phi(\Gamma)}r_*M$ is finitely generated over $gr^\Phi ZA_n(K)$. Let us show the latter.

Consider the finite exhaustive filtration of $gr^{\Phi(\Gamma)}r_*M$ by graded $gr^\Phi ZA_n(K)$-submodules, $$0 =(gr^{\Phi(\Gamma)}r_*M)_0 \subset (gr^{\Phi(\Gamma)}r_*M)_1 \subset...\subset (gr^{\Phi(\Gamma)}r_*M)_p= gr^{\Phi(\Gamma)}r_*M.$$ It is defined as follows, $\forall 0\leq i\leq p, \forall l\in \mathbb{Z},$ let $pm$ the greatest integer multiple of $p$ such that $pm\leq l.$ Then $(gr^{\Phi(\Gamma)}r_*M)_i \cap \Phi(\Gamma)_lr_*M/\Phi(\Gamma)_{l-1}r_*M$ is the image of the map $\Gamma_{p(m-1)+i}M 
\to \Phi(\Gamma)_lr_*M/\Phi(\Gamma)_{l-1}r_*M$ induced by the inclusion $\Gamma_{p(m-1)+i}M \subset \Gamma_{pm}M=:\Phi(\Gamma)_lr_*M.$ 
Let  $gr(gr^{\Phi(\Gamma)}r_*M):= \bigoplus\limits_{i=1}^{i=p}(gr^{\Phi(\Gamma)}r_*M)_i/(gr^{\Phi(\Gamma)}r_*M)_{i-1}$ be the associated graded $gr^\Phi ZA_n(K)$-module.  
Note furthermore that the module $gr^{\Gamma}M$ seen as a module over $gr^\Phi ZA_n(K) \hookrightarrow gr^{\mathcal{B}}A_n(K)$ decomposes into a direct sum of graded $gr^\Phi ZA_n(K)$-submodules $\bigoplus\limits_{i=1}^{i=p}(gr^{\Gamma}M)_i,$ where $(gr^{\Gamma}M)_i:=\bigoplus_{l\in\mathbb{Z}}gr_{pl+i}^{\Gamma}M.$ Let $F_*gr^{\Gamma}M$ be the graded $gr^\Phi ZA_n(K)$-module $\bigoplus\limits_{i=1}^{i=p}(gr^{\Gamma}M)_i[i-p],$ where $[\bullet]$ denotes the degree shift.
Then, $\forall 1\leq i\leq p,$ there is an isomorphism of graded $gr^\Phi ZA_n(K)$-modules $$(gr^{\Phi(\Gamma)}r_*M)_i/(gr^{\Phi(\Gamma)}r_*M)_{i-1}\to (gr^{\Gamma}M)_i[i-p].$$ Indeed in degree $l\in \mathbb{Z},$ with $pm$ the greatest integer multiple of $p$ such that $pm\leq l,$ it is induced by the quotient map $\Gamma_{p(m-1)+i}M \to gr_{p(m-1)+i}^{\Gamma}M,$ the latter being equal to $((gr^{\Gamma}M)_i)_{p(m-1)+i} = ((gr^{\Gamma}M)_i[i-p])_{pm},$ where the outermost index refers to the homogeneous component of a graded $gr^\Phi ZA_n(K)$-module.  
These finally assemble into an isomorphism of graded $gr^\Phi ZA_n(K)$-modules 
\begin{equation} \label{eq: grZ}
gr(gr^{\Phi(\Gamma)}r_*M) \simeq F_*gr^{\Gamma}M.
\end{equation} 

We conclude by noting that the $gr^\Phi ZA_n(K)$-module $F_*gr^{\Gamma}M$ is finitely generated. Indeed, by hypothesis, the $gr^{\mathcal{B}}A_n(K)$-module $gr^{\Gamma}M$ is finitely generated. Hence it is finitely generated as a $gr^\Phi ZA_n(K)$-module since $gr^{\mathcal{B}}A_n(K)\cong K[x_1,...,x_n,y_1,...,y_n]$ is a finitely generated module over $K[x_1^p,...,x_n^p,y_1^p,...,y_n^p] \cong gr^\Phi ZA_n(K).$ Thus $F_*gr^{\Gamma}M$ is a finitely generated $gr^\Phi ZA_n(K)$-module as direct summands of a finitely generated module are finitely generated. So, by the above isomorphism, the $gr^\Phi ZA_n(K)$-module $gr(gr^{\Phi(\Gamma)}r_*M)$ is finitely generated. Consequently the finite exhaustive filtration of $gr^{\Phi(\Gamma)}r_*M$ has finitely generated subquotients and hence $gr^{\Phi(\Gamma)}r_*M$ is a finitely generated $gr^\Phi ZA_n(K)$-module. This finishes the proof of the lemma.\end{proof}

Let $M$ be a left $A_n(K)$-module and let $\Gamma$ be a good filtration of $M.$ By the above lemma \ref{lmm:pGamma is good}, the Rees module of $(r_*M,p\Gamma)$ is a finitely generated graded module over the Rees ring $R_n(\mathcal{C})\simeq K[t_0,x_1^p,...,x_n^p,\partial_1^p,...,\partial_n^p].$ Thus it has a Hilbert polynomial $\mathcal{H}_{\mathcal{R}(r_*M,p\Gamma)}.$ In the following proposition, we express $\mathcal{H}_{\mathcal{R}(r_*M,p\Gamma)}$ in terms of the Hilbert polynomial
$\mathcal{H}_{M,\Gamma}$ of $(M,\Gamma).$ 

\begin{prop} \label{prop:Hilbert polynomial of pGamma}

Let $M$ be a left $A_n(K)$-module and let $\Gamma$ be a good filtration of $M.$ The Hilbert polynomial $\mathcal{H}_{\mathcal{R}(r_*M,p\Gamma)}(t)$ of the Rees module of $(r_*M,p\Gamma)$ is $\mathcal{H}_{M,\Gamma}(pt).$ In particular, the degree of $\mathcal{H}_{\mathcal{R}(r_*M,p\Gamma)}$ is the dimension $d(M)$ of $M$ and its leading coefficient times $d(M)!$ is $e(M)p^{d(M)},$ where $e(M)$ is the multiplicity of $M.$

\end{prop}
\begin{proof} For $l$ large enough, on the one hand the Hilbert polynomial $\mathcal{H}_{M,\Gamma}(l)$ coincides with the function $l\mapsto \dim_K\Gamma_lM$ and on the other $\mathcal{H}_{\mathcal{R}(r_*M,p\Gamma)}(l)$ coincides with $l\mapsto \dim_K(p\Gamma)_lM=\dim_K\Gamma_{pl}M.$ The proposition immediately follows.\end{proof}

\subsection{Conclusion} \label{subsec:bound} Here we obtain the estimates mentioned at the beginning of the section.  

First, we would like to recall a well-known geometric interpretation of the Rees ring and module. 
We use the notations of \cite[\S2]{EGAII} for projective schemes.
\begin{prop} \label{prop: Rees}
Using the notations of Lemma \ref{lmm: Rees ring}, we have $$D_+(t_0) \xhookrightarrow{j} Proj(R_n(\mathcal{C})) \xhookleftarrow{i} V_+(t_0)$$ where $j$ is an open embedding and $i$ is its closed complement, and

\begin{enumerate}[(a)]
\item \label{item: a} $D_+(t_0)\cong \Spec(ZA_n(K))$
\item \label{item: b}$V_+(t_0)\cong Proj(gr^{\mathcal{C}}ZA_n(K))$
\end{enumerate}

\end{prop}

\begin{proof} (\ref{item: a}) By definition, $D_+(t_0)= \Spec(R_n(\mathcal{C})_{(t_0)}).$ Hence (\ref{item: a}) follows immediately from (\ref{item: sum}) of Lemma \ref{lmm: Rees ring}.

(\ref{item: b}) By definition, $V_+(t_0)= Proj(R_n(\mathcal{C})/t_0R_n(\mathcal{C})).$ Hence (\ref{item: b}) follows from (\ref{item: fiber}) of Lemma \ref{lmm: Rees ring}.
\end{proof}

Making these identifications, let $G$ be a good filtration on a finitely generated $(ZA_n(K),\mathcal{C})$-module $M.$ Then the coherent sheaf $\widetilde{\mathcal{R}(M,G)}$ on $Proj(R_n(\mathcal{C}))$ extends the coherent sheaf $\widetilde{M}$ on $\Spec(ZA_n(K))$ and its restriction to the complement $Proj(gr^{\mathcal{C}}ZA_n(K))$ of $\Spec(ZA_n(K))$ is isomorphic to $\widetilde{gr^GM}$. Finally, one easily sees that the support of $\widetilde{\mathcal{R}(M,G)}$ is the closure of $supp(\widetilde{M})$ in $Proj(R_n(\mathcal{C}))$.
Note that here \label{page:Rees} $$Proj(R_n(\mathcal{C})) \simeq Proj(K[t_0,x_1^p,...,x_n^p,\partial_1^p,...,\partial_n^p])\simeq\mathbb{P}^{2n}_K, $$ $$\Spec(ZA_n(K)) \simeq \Spec(K[x_1^p,...,x_n^p,\partial_1^p,...,\partial_n^p])\simeq \mathbb{A}^{2n}_K $$ $$Proj(gr^{\mathcal{C}}ZA_n(K))\simeq Proj(K[x_1^p,...,x_n^p,\partial_1^p,...,\partial_n^p])\simeq\mathbb{P}^{2n-1}_K$$

The leading coefficient of the Hilbert polynomial of $\widetilde{\mathcal{R}(M,G)}$ is related to the top-dimensional irreducible components of its support through the following:

\begin{prop} \label{prop:mu, ranks and degrees}

Let $Y \overset{i}\hookrightarrow \mathbb{P}^{m}_K$ be a closed subscheme and let $\mathcal{F}$ be a coherent sheaf of dimension $d$ on $Y$. Let the degree of $\mathcal{F}$ with respect to $i$ be $\mu(\mathcal{F}):=d!a_d$ where $a_d$ is the leading coefficient of the Hilbert polynomial of $\mathcal{F}$ with respect to $i$. Then $$\Sigma_z \mathrm{rk}_z(\mathcal{F})\deg(\overline{\{z\}}) \leq\mu(\mathcal{F})$$ where the sum is over the generic points of the $d$-dimensional irreducible components of $supp(\mathcal{F})$, $\mathrm{rk}_z(\mathcal{F}):=\dim_{k(z)}(\mathcal{F}_z\otimes k(z))$ and $\deg(\overline{\{z\}})$ is the degree of $\overline{\{z\}}^{red}$ with respect to $i$.

\end{prop}
\begin{proof} By \cite[Lemma B.4]{Kleiman} and \cite[Proposition 5.3.1]{EGAIV}, $\mu(\mathcal{F})=\Sigma_z \lg_{\mathcal{O}_{Y,z}}(\mathcal{F}_z)\mu(\mathcal{O}_{\overline{\{z\}}^{red}}),$ where $\lg$ denotes the length, 
summing over the generic points of the $d$-dimensional irreducible components of $supp(\mathcal{F}).$ Let $z$ be as above, then by additivity of the length under short exact sequences $\lg_{\mathcal{O}_{Y,z}}(\mathcal{F}_z)\geq \lg_{k(z)}(\mathcal{F}_z\otimes k(z))=\dim_{k(z)}(\mathcal{F}_z\otimes k(z))=:\mathrm{rk}_z(\mathcal{F}).$ The proposition follows as $\deg(\overline{\{z\}}):=\mu(\mathcal{O}_{\overline{\{z\}}^{red}}).$\end{proof}

\begin{thm} \label{thm:bound}

Let $S$ be an integral scheme dominant and of finite type over $\mathbb{Z}$ and let $M$ be a coherent left $D_{\mathbb{A}^n_S/S}$-module. Let $\mu$ be the generic point of $S.$ Suppose that $M_{\mu}$ is a nonzero holonomic left $D_{\mathbb{A}^n_{k(\mu)}}$-module. Then there is a dense open subset $U$ of $S$ such that for each closed point $u\in U$ and each $z$ generic point of an irreducible component of $p$-supp($M_u$) $$\deg(\overline{\{z\}})\leq e(M_\mu) \text{ and }\mathrm{rk}_z(M_u) \leq e(M_\mu)p^n$$ where $e(M_\mu)$ is the multiplicity for the Bernstein filtration of $M_{\mu},$ $\deg(\overline{\{z\}})$ is the degree of the reduced closure of the image of $z$ in $\mathbb{P}^{2n}_{k(u)}$ by the open immersion of the Rees construction and $\mathrm{rk}_z(M_u):=\dim_{k(z)}((F_*{M_u})_z\otimes k(z)).$

\end{thm}
\begin{proof} The proof reduces to the case of an integral and affine $S = \Spec(R)$. We may thus consider that $M$ is a finitely generated left $A_n(R)$-module. By Lemma \ref{lmm:fibers dimension and multiplicity}, there is a dense open subset $U_e$ of $S$ such that for each closed point $u\in U_e$, $d(M_u)=n$ and $e(M_u)=e(M_\mu).$ 

Let $u\in U_e$ with $p = char(k(u))$ and let $\Gamma$ be a good filtration on the left $A_n(k(u))$-module $M_u$. By Proposition \ref{prop:Hilbert polynomial of pGamma}, $\widetilde{\mathcal{R}(M_u,p\Gamma)}$ is of dimension $n$ and $\mu(\widetilde{\mathcal{R}(M_u,p\Gamma)})=e(M_\mu)p^n$. Moreover $supp(\widetilde{\mathcal{R}(M_u,p\Gamma)})$ is the closure $\overline{p \text{-}supp(M_u)}$ of $p$-$supp(M_u),$ in which $p$-$supp(M_u)=\overline{p \text{-}supp(M_u)}\cap \Spec(ZA_n(K))$ is open. Hence Proposition \ref{prop:mu, ranks and degrees} implies that $\Sigma_z \mathrm{rk}_z(M_u)\deg(\overline{\{z\}})\leq e(M_\mu)p^n,$ where the sum is over the generic points of the $n$-dimensional irreducible components of $p$-supp($M_u$). 

By the equidimensionality of the $p$-supports (Theorem \ref{thm:purity}), there is a dense open subset $U\subset U_e$ such that for all closed points $u$ of $U,$ all the irreducible components of $p$-$supp(M_u)$ are of dimension $n.$ Hence we deduce that, for each $z$ generic point of an irreducible component of $p$-supp($M_u$), $\mathrm{rk}_z(M_u)\deg(\overline{\{z\}})\leq e(M_\mu)p^n$. This implies the second estimate of the theorem.

Finally, let $u$ be a closed point in $U$ of characteristic $p.$ By \cite[Theorem 2.2.3]{BMR}, $F_*M_u$ is a left module over an Azumaya algebra of rank $p^{2n}.$ Hence $(F_*M_u)_z\otimes \overline{k(z)}$ is by \cite[Th\'eor\`eme 5.1 (i)]{BrauerI} a left module over the algebra of $p^n\times p^n$ matrices $M_{p^n}(\overline{k(z)}),$ where $\overline{k(z)}$ is an algebraic closure of $k(z).$ Therefore by Morita theory there is a finite dimensional $\overline{k(z)}$-vector space $V$ such that $(F_*M_u)_z\otimes \overline{k(z)} \simeq\overline{k(z)}^{p^n} \otimes_{\overline{k(z)}}V,$ where $\overline{k(z)}^{p^n}$ is the standard left $M_{p^n}(\overline{k(z)})$-module. In particular $\mathrm{rk}_z(M_u):=\dim_{k(z)}((F_*M_u)_z\otimes k(z))=\dim_{\overline{k(z)}}((F_*M_u))_z\otimes \overline{k(z)})$ is divisible by $p^n.$ The first estimate of the theorem follows.\end{proof}

\begin{rem} The first estimate was conjectured in \cite[Conjecture 1]{Kontsevich}.
\end{rem}

\section{The Brauer group and differential forms} \label{sec:Brauer group and forms}

Here we prove, in a first part, that the Azumaya algebra of differential operators splits on the regular locus of the $p$-support of a holonomic $\mathcal{D}$-module, for $p$ large enough. See Theorem \ref{thm:splitting}.

In a second part, we recollect some facts of differential calculus in positive characteristic. In particular, we consider the $p$-curvature exact sequence  and a map arising from it  
which sends 1-forms to the Brauer group. The image of the canonical form is the class of the Azumaya algebra of differential operators.  
In view of the first part of the section and for later use in the proof of our main theorem, we describe its kernel. 

\subsection{Splittings of Azumaya algebras on the support of their modules} \label{subsec: splitting}
Let $Y$ be a scheme and let $\mathcal{A}$ be an Azumaya algebra on $Y$. Since $\mathcal{A}$ is a coherent $\mathcal{O}_Y$-module, it is a coherent Noetherian ring and a left $\mathcal{A}$-module is coherent if and only if it is coherent as an $\mathcal{O}_Y$-module. Recall that an Azumaya algebra $\mathcal{A}$ is said to split on $Y$ if its class $[\mathcal{A}]$ in the Brauer group $Br(Y)$ of $Y$ is trivial.  

Let $M$ be a coherent left $\mathcal{A}$-module and let $z$ be the generic point of an irreducible component of the support of the coherent $\mathcal{O}_Y$-module $M.$ The next proposition shows that the rank $\mathrm{rk}_z(M)$ of $M$ at $z$ constrains  
the order of $[\mathcal{A}|_{(\overline{\{z\}}^{red})^{reg}}]$ in $Br((\overline{\{z\}}^{red})^{reg}).$ 

We first prove a lemma.

\begin{lmm} \label{lmm:Azumaya acting on a vector bundle}

Let $Y$ be a scheme and let $\mathcal{A}$ be an Azumaya algebra of rank $r^2$ on $Y$. Suppose that $\mathcal{A}$ acts on the left on a locally free sheaf $\mathcal{V}$ of rank $v$. Then $r$ divides $v=lr$ and $l[\mathcal{A}]=0$ in $Br(Y).$

\end{lmm}
\begin{proof} By hypothesis there is a morphism of $\mathcal{O}_Y$-algebras $\mathcal{A} \to \mathcal{E}nd_{\mathcal{O}_Y}(\mathcal{V}),$ sending $1$ to $1.$ It is injective by \cite[Proposition 0.5.5.4]{EGAI} since the fiber of $\mathcal{A}$ at each point of $Y$ is a simple algebra by \cite[Th\'eor\`eme 5.1 (i)]{BrauerI}. Therefore one may view $\mathcal{A}$ as a subalgebra of $\mathcal{E}nd_{\mathcal{O}_Y}(\mathcal{V})$ and in particular consider the commutant $\mathcal{C}_{\mathcal{E}nd_{\mathcal{O}_Y}(\mathcal{V})}(\mathcal{A})$  
of $\mathcal{A}$ in $\mathcal{E}nd_{\mathcal{O}_Y}(\mathcal{V}),$ which is a coherent subalgebra of $\mathcal{E}nd_{\mathcal{O}_Y}(\mathcal{V}).$ By \cite[Theorem 3.3]{AusGold}, the natural morphism of $\mathcal{O}_Y$-algebras $\mathcal{A}\otimes_{\mathcal{O}_Y}\mathcal{C}_{\mathcal{E}nd_{\mathcal{O}_Y}(\mathcal{V})}(\mathcal{A})\to \mathcal{E}nd_{\mathcal{O}_Y}(\mathcal{V})$ is an isomorphism and $\mathcal{C}_{\mathcal{E}nd_{\mathcal{O}_Y}(\mathcal{V})}(\mathcal{A})$ is an Azumaya algebra on $Y$. Hence by the behaviour of ranks under tensor products, $\mathcal{C}_{\mathcal{E}nd_{\mathcal{O}_Y}(\mathcal{V})}(\mathcal{A})$ is of constant rank $l^2,$ such that $v=lr$. By definition of the Brauer group $0=[\mathcal{E}nd_{\mathcal{O}_Y}(\mathcal{V})]=[\mathcal{A}]+[\mathcal{C}_{\mathcal{E}nd_{\mathcal{O}_Y}(\mathcal{V})}(\mathcal{A})]$ in $Br(Y).$ The lemma follows since for each Azumaya algebra $\mathcal{B}$ of rank $n^2$ on $Y$, $n[\mathcal{B}]=0$ in $Br(Y)$ (\cite[\S 2]{BrauerI}) giving $0=l[\mathcal{A}]+l[\mathcal{C}_{\mathcal{E}nd_{\mathcal{O}_Y}(\mathcal{V})}(\mathcal{A})]=l[\mathcal{A}].$ 
\end{proof}

\begin{prop} \label{prop:order on the support}

Suppose that $Y$ is of finite type over a field $K$. Let $\mathcal{A}$ be an Azumaya algebra of rank $r^2$ on $Y$, let $M$ be a coherent left $\mathcal{A}$-module and let $z$ be the generic point of an irreducible component $\overline{\{z\}}$ of $supp(M)$. Then $r$ divides $\mathrm{rk}_z(M)=l_z(M)r$ and $$l_z(M)[\mathcal{A}|_{(\overline{\{z\}}^{red})^{reg}}]=0$$ in $Br((\overline{\{z\}}^{red})^{reg}).$

\end{prop}
\begin{proof} Since the vector space $M_z\otimes k(z)$ is of dimension $\mathrm{rk}_z(M)$ and acted upon on the left by the rank $r^2$ Azumaya algebra $\mathcal{A}_z\otimes k(z),$ Lemma \ref{lmm:Azumaya acting on a vector bundle} implies that $\mathrm{rk}_z(M)=l_z(M)r$ and $l_z(M)[\mathcal{A}_z\otimes k(z)]=0$ in $Br(k(z)).$  
Moreover the restriction of $\mathcal{A}$ to the regular locus of the irreducible component satisfies $\mathcal{A}_z\otimes k(z)\cong (\mathcal{A}|_{\overline{\{z\}}^{red}})_z\otimes k(z) \cong (\mathcal{A}|_{(\overline{\{z\}}^{red})^{reg}})_z\otimes k(z).$ The proposition then follows from the canonical embedding $Br((\overline{\{z\}}^{red})^{reg}) \hookrightarrow Br(k(z)),$ see \cite[IV Corollary 2.6]{Milne}.\end{proof}

The following combines Proposition \ref{prop:order on the support} with the second estimate of Theorem \ref{thm:bound}.   

\begin{thm} \label{thm:splitting in affine space case}

Let S be an integral scheme dominant and of finite type over $\mathbb{Z}$ and let $M$ be a coherent left $D_{\mathbb{A}^n_S/S}$-module. Let $\mu$ be the generic point of $S.$ Suppose that $M_{\mu}$ is a nonzero holonomic left $D_{\mathbb{A}^n_{k(\mu)}}$-module. Then there is a dense open subset $U$ of $S$ such that for each closed point $u\in U$ and each $z$ generic point of an irreducible component of $p$-supp($M_u$), the Azumaya algebra $F_*D_{\mathbb{A}^n_{k(u)}}$ on $T^*_{(\mathbb{A}^n_{k(u)})'}$ splits on $(\overline{\{z\}}^{red})^{reg}.$

\end{thm}
\begin{proof} By Theorem \ref{thm:bound} and using its notations, there is a dense open subset $U_b$ of $S$ such that for each closed point $u\in U_b$ and each $z$ generic point of an irreducible component of $p$-supp($M_u$), $\mathrm{rk}_z(M_u) \leq e(M_\mu)p^n$ where $p$ is the characteristic of the residue field $k(u).$ Thus using Proposition \ref{prop:order on the support} and its notations, we have   
$l_z(M_u)\leq e(M_\mu)$ and $l_z(M_u)[F_*D_{\mathbb{A}^n_{k(u)}}|_{(\overline{\{z\}}^{red})^{reg}}]=0$ in $Br((\overline{\{z\}}^{red})^{reg}).$ Note that by definition $l_z(M_u)\neq 0$ and hence for $u \in U$ the open dense subset of $U_b$ defined by inverting all the primes $\leq e(M_\mu)$, $l_z(M_u)$ and $p^n$ are coprime, that is there are integers $a$ and $b$ such that $1=al_z(M_u)+bp^n.$ Since $F_*D_{\mathbb{A}^n_{k(u)}}$ is of rank $p^{2n}$, $p^n[F_*D_{\mathbb{A}^n_{k(u)}}|_{(\overline{\{z\}}^{red})^{reg}}]=0$ by \cite[\S 2]{BrauerI} and the theorem follows from $[F_*D_{\mathbb{A}^n_{k(u)}}|_{(\overline{\{z\}}^{red})^{reg}}]
=1[F_*D_{\mathbb{A}^n_{k(u)}}|_{(\overline{\{z\}}^{red})^{reg}}]
$

$=(al_z(M_u)+bp^n)[F_*D_{\mathbb{A}^n_{k(u)}}|_{(\overline{\{z\}}^{red})^{reg}}]
$

$=al_z(M_u)[F_*D_{\mathbb{A}^n_{k(u)}}|_{(\overline{\{z\}}^{red})^{reg}}]+bp^n[F_*D_{\mathbb{A}^n_{k(u)}}|_{(\overline{\{z\}}^{red})^{reg}}]$

$=0$ in $Br((\overline{\{z\}}^{red})^{reg}).$\end{proof}

As a corollary, we have the following.
  
\begin{thm} \label{thm:splitting}
Let $S$ be an integral scheme dominant and of finite type over $\mathbb{Z}$, let $X$ be a smooth $S$-scheme of relative dimension $n$ and let $M$ be a coherent left $D_{X/S}$-module. Suppose that the fiber of $M$ at the generic point $\mu$ of $S$ is a holonomic left $D_{X_\mu}$-module. Then there is a dense open subset $U$ of $S$ such that for each closed point $u$ of $U,$ the Azumaya algebra $F_*D_{X_u}$ on $T^*_{X'_u}$ splits on the regular locus of the $p$-support $p\text{-}supp(M_u)^{reg}.$
\end{thm}
\begin{proof} If $M_\mu = 0,$ the theorem is trivial by  
Remark \ref{rem:generic fiber zero}. We thus suppose that $M_\mu$ does not vanish. 

By the canonical injection of the Brauer group for the Zariski site (i.e. the classes of Azumaya algebras which are Zariski-locally isomorphic to an algebra of matrices) into the Zariski cohomology $H^2(Y,\mathcal{O}_Y^\times)$ (\cite[(2.1)]{BrauerI}),  
the case $i=2$ of Lemma \ref{lmm:vanishing of Zariski cohomology} below implies that on a regular (Noetherian) scheme for an Azumaya algebra to be split is a Zariski-local condition. 

Therefore by Remark \ref{rem:reduction to affine}, one may further assume that $S$ and $X$ are regular integral and affine and in particular that there is a closed immersion $X \overset{f}\hookrightarrow \mathbb{A}^m_S$ over $S.$ 

Let us specialize to a closed point $u$ of positive characteristic $p$ of $S.$ Since $f'_d$ is smooth by Lemma \ref{lmm:cotangent of an immersion}, the description of $p\text{-}supp(H^0(f_u)_+M_u)$ given in the proof of Proposition \ref{prop:reduction to affine space}, of which we use the notations,  
implies that $p\textit{-supp}(H^0(f_u)_+M_u)^{reg}= ({f_u'})_{\pi}\circ({f_u'})_d^{-1}(p\text{-}supp(M_u)^{reg}).$ 
 
By \cite[Proposition 3.7]{BezrBr}, ${f'_d}^*(F_*D_{X_u})$ splits on ${f'_d}^{-1}(p\text{-}supp(M_u)^{reg})$ if $F_*D_{\mathbb{A}^m_{k(u)}}$ splits on the regular locus of
$p\textit{-supp}(H^0(f_u)_+M_u).$ Moreover the pullback of Brauer classes ${f'_d}^*$  
induces an injective morphism of Brauer groups $Br(p\text{-}supp(M_u)^{reg}) \hookrightarrow Br({f'_d}^{-1}(p\text{-}supp(M_u)^{reg})).$ Indeed $f'_d$ Zariski-locally admits a section by Lemma \ref{lmm:cotangent of an immersion} and, on the regular locus, being split is Zariski-local as we explained above.  
So $F_*D_{X_u}$ splits on the regular locus of $p\text{-}supp(M_u)$ if $F_*D_{\mathbb{A}^m_{k(u)}}$ splits on $p\textit{-supp}(H^0(f_u)_+M_u)^{reg}.$ 

Finally,  
an Azumaya algebra splits on a regular Noetherian scheme if and only if it splits on its irreducible components. Thus  
the theorem follows from Theorem \ref{thm:splitting in affine space case}.\end{proof}

\begin{lmm} \label{lmm:vanishing of Zariski cohomology}

Let $Y$ be a noetherian scheme. If $Y$ is locally factorial, i.e. its local rings are UFD, then the Zariski cohomology $H^i(Y,\mathcal{O}_Y^\times)=0,$ for all $i\geq 2.$ 

\end{lmm}
\begin{proof} By definition of the sheaf $\mathcal{D}iv_Y$ of Cartier divisors there is an exact sequence of abelian sheaves $0\to\mathcal{O}_Y^\times\to\mathcal{K}_Y^\times\to\mathcal{D}iv_Y\to 0$ on $Y,$ where $\mathcal{K}_Y$ is the sheaf of meromorphic functions and $\mathcal{O}_Y^\times\to\mathcal{K}_Y^\times$ is the natural injection. If $Y$ is locally factorial then it is the sum of its (finitely many) irreducible components, each of which is integral, by \cite[Proposition 4.5.5]{EGAI}. Hence if $Y_i \overset{f_i}\hookrightarrow Y$ is the open immersion of the $i$-th irreducible component, then $\mathcal{K}_Y^\times \cong \Pi_i {f_i}_*\mathcal{K}_{Y_i}^\times.$ Moreover $\mathcal{K}_{Y_i}^\times$ is isomorphic to the constant sheaf associated to $k(y_i)^\times$ for $y_i$ the generic point of $Y_i.$ In particular $\mathcal{K}_Y^\times$ is flasque. Since $\mathcal{D}iv_Y$ is flasque by \cite[Corollaire 21.6.11]{EGAIV}, $\mathcal{K}_Y^\times\to\mathcal{D}iv_Y$ is a flasque right resolution of $\mathcal{O}_Y^\times,$ vanishing in degrees $\geq2.$ The result follows as sheaf cohomology may be computed using flasque resolutions.\end{proof}

\subsection{The Brauer group via the $p$-curvature exact sequence} \label{subsec:Brauer group via the p-curvature sequence}
Let $Y$ be a smooth  
scheme over a perfect field $K$ of positive characteristic $p.$ Let $Y \xrightarrow{F} Y'$ be the relative Frobenius morphism and let $Y'\xrightarrow{W}Y$ be the projection. 
We denote by $\bigoplus_{i\in \mathbb{Z}} \Omega^{i}_{Y'} \xrightarrow{C_Y^{-1}} \bigoplus_{i\in \mathbb{Z}} \mathcal{H}^{i}(F_*\Omega^{\bullet}_Y)$ the Cartier isomorphism, see \cite[Theorem 7.2]{Katz}.

\begin{defi} The Cartier operator is the composed morphism $$C_Y: \bigoplus_{i\in \mathbb{Z}} Z^i(F_*\Omega^{\bullet}_Y) \xrightarrow{\pi} \bigoplus_{i\in \mathbb{Z}} \mathcal{H}^i(F_*\Omega^\bullet_Y) \xrightarrow{\text{inverse of }C_Y^{-1}} \bigoplus_{i\in \mathbb{Z}} \Omega^i_{Y'}, $$ where $Z^i(F_*\Omega^\bullet_Y):= \ker (F_*\Omega^i_Y \xrightarrow{F_*d} F_*\Omega^{i+1}_Y)$ is the sheaf of closed $i$-forms and $\pi$ is the quotient morphism. 
\end{defi}

Recall that there is an exact sequence of \'{e}tale sheaves on $Y':$ 

\begin{equation} \label{equ: p-curvature}
0\to \mathbb{G}_{\,m/Y'} \xrightarrow{F^*} F_*\mathbb{G}_{\,m/Y} \xrightarrow{F_*dlog} F_*Z^1(\Omega^{\bullet}_Y) \xrightarrow{W^{\star}-C_Y} \Omega^{1}_{Y'} \to 0
\end{equation} 
where $dlog(y):= dy/y$ for all local sections $y$ of $\mathbb{G}_{\,m/Y}$
and $W^{\star}$ is induced by the pullback on forms, see \cite[Corollaire 0.2.1.18]{Illusie dR-W}.

\begin{defi} \label{defi: p-curvature}
The exact sequence (\ref{equ: p-curvature}) is the $p$-curvature exact sequence and the morphism $F_*Z^1(\Omega^{\bullet}_Y) \xrightarrow{W^{\star}-C_Y} \Omega^{1}_{Y'}$ is the $p$-curvature operator.
\end{defi}

The $p$-curvature exact sequence (\ref{equ: p-curvature}) decomposes into two short exact sequences of \'{e}tale sheaves on $Y':$

\begin{equation}
0\to \mathbb{G}_{\,m/Y'} \xrightarrow{F^\ast} F_*\mathbb{G}_{\,m/Y} \xrightarrow{F_*dlog} \text{Im}F_*dlog\to 0
\end{equation}
\begin{equation}
0\to \coker F^*\xrightarrow{F_*dlog} F_*Z^1(\Omega^{\bullet}_Y) \xrightarrow{W^{\star}-C_{Y}} \Omega^{1}_{Y'} \to 0
\end{equation}

Composing the coboundary morphisms of the corresponding \'{e}tale cohomology long exact sequences, one deduces a morphism:
$$\Psi_Y: H^0(Y',\Omega^{1}_{Y'}) \to H^1(Y',\coker{F}^*\cong \text{Im}F_*dlog) \to H^2(Y',\mathbb{G}_{\,m/Y'})$$ 

\begin{prop} \label{prop: phi}
The morphism $H^0(Y',\Omega^{1}_{Y'}) \xrightarrow{\Psi_Y} H^2(Y',\mathbb{G}_{\,m/Y'})$ factors uniquely through the canonical embedding $Br(Y')_p\hookrightarrow H^2(Y',\mathbb{G}_{\,m/Y'}),$ where $Br(Y')_p$ is the kernel of multiplication by $p$ in $Br(Y').$\end{prop}

\begin{proof} This is well-known. By construction $\Psi_Y$ factors through
$$H^0(Y',\Omega^{1}_{Y'}) \twoheadrightarrow \coker H^0(W^{\star}-C_{Y}) \to \ker H^2(F^*)\hookrightarrow H^2(Y',\mathbb{G}_{\,m/Y'}),$$ the end maps being the natural ones. Moreover by \cite[Proposition 2.1]{Hoobler}, the image of $\ker H^2(F^*)\hookrightarrow H^2(Y',\mathbb{G}_{\,m/Y'})$ is the kernel $H^2(Y',\mathbb{G}_{\,m/Y'})_p$ of multiplication by $p,$ the latter being the image of the canonical embedding $Br(Y')_p\hookrightarrow H^2(Y',\mathbb{G}_{\,m/Y'}).$\end{proof}

\begin{defi} \label{defi: phi}The morphism deduced from Proposition \ref{prop: phi} is denoted $H^0(Y',\Omega^{1}_{Y'}) \xrightarrow{\phi_Y} Br(Y')_p.$ 
\end{defi}
 
\begin{rem}\label{rem:other description of phi}

Here is another description of $\phi_Y,$ from \cite[Remark 4.3]{Ogus-Vologodsky}: 

Let $\alpha \in H^0(Y',\Omega^{1}_{Y'}),$ then $\phi_Y(\alpha)=[s_{\alpha}^*(F_*D_{Y})] \in Br(Y'),$ where $Y' \overset{s_{\alpha}}\to T^*_{Y'}$ is the section of $T^*_{Y'}/Y'$ corresponding to $\alpha.$ 

\end{rem}

The morphism $\phi_Y$ depends functorially on $Y:$

\begin{lmm} \label{lmm:functoriality of phi}

Let $Z\xrightarrow{f} Y$ be a $K$-morphism of smooth $K$-schemes, $Z'\xrightarrow{f'} Y'$ its base-change by Frobenius and let $\alpha \in H^0(Y',\Omega^1_{Y'}).$ Then $f'^*\phi_Y(\alpha)=\phi_Z((f'^*)^{ad}\alpha),$ where $f'^*$ on the left (resp. on the right) is the pullback of classes in the Brauer group (resp. pullback of forms) by $f'$.

\end{lmm}
\begin{proof} Set $\mathcal{D}_Y:= F_*D_{Y}$ and $\mathcal{D}_Z:= F_*D_{Z}.$ With the notations of the Appendix, we have:
By Remark \ref{rem:other description of phi}, $\phi_Z((f'^*)^{ad}\alpha)=[s_{(f'^*)^{ad}\alpha}^*(\mathcal{D}_Z)].$
And $[s_{(f'^*)^{ad}\alpha}^*(\mathcal{D}_Z)]=[(f'_d \circ (Z\times_Ys_{\alpha}))^*(\mathcal{D}_Z)],$ since $s_{(f'^*)^{ad}\alpha}=f'_d \circ Z\times_Ys_{\alpha},$ by Remark \ref{rem: pullback section}. Moreover by \cite[Proposition 3.7]{BezrBr}, $[{f'_d}^*(\mathcal{D}_Z)]=[{f'_{\pi}}^*(\mathcal{D}_Y)].$ Hence $[(f'_d \circ (Z\times_Ys_{\alpha}))^*(\mathcal{D}_Z)]=[(Z\times_Ys_{\alpha})^*{f'_d}^*(\mathcal{D}_Z)]
=[(Z\times_Ys_{\alpha})^*{f'_{\pi}}^*(\mathcal{D}_Y)]
=[(f'_{\pi} \circ (Z\times_Ys_{\alpha}))^*(\mathcal{D}_Y)]=[(s_{\alpha}\circ f')^*(\mathcal{D}_Y)],$  
using the equality $f'_{\pi} \circ (Z\times_Ys_{\alpha}) = s_{\alpha}\circ f'.$ Finally by Remark \ref{rem:other description of phi}, the last term is equal to $f'^*\phi_Y(\alpha),$ as stated.\end{proof}
 
By construction $\phi_Y$ factors through $\coker H^0(W^{\star}-C_{Y})$ and we denote the resulting map by $\coker H^0(W^{\star}-C_{Y}) \xrightarrow{\overline{\phi_Y}} Br(Y')_p.$ 
The following proposition provides information on the kernel of $\phi_Y.$

\begin{prop} \label{prop:kernel of phi}

Suppose further that $Y$ is affine. Then there is an exact sequence, commuting with restriction to affine open subsets:
$$Pic(Y)\to \coker H^0(W^{\star}-C_{Y}) \xrightarrow{\overline{\phi_Y}} Br(Y')_p \to 0.$$

\end{prop}
\begin{proof} This is a special case of \cite[Corollary 1.7]{Hoobler}.\end{proof}

Let us now specialize to the case $Y=T^*_X,$ with $X$ is a smooth  
$K$-scheme. Let $X \xrightarrow{F_X} X'$ be the Frobenius morphism and $X' \xrightarrow{W_X} X$ the projection.

\begin{rem} \label{rem: T^*X'} The pullback of forms $W_X^*\Omega^1_X \xrightarrow{W_X^*}\Omega^1_{X'}$ is an isomorphism and induces an isomorphism ${T^*_X}'\tilde{\to} T^*_{X'}.$ We use this isomorphism  
to identify ${T^*_X}'$ and $T^*_{X'}.$
\end{rem}

We have the following description of the Brauer class of the algebra of differential operators in terms of $\phi:$

\begin{prop} \label{prop:canonical form maps to algebra of differential operators}
Let $\theta_{X'}$ be the canonical $1$-form on $T^*_{X'},$ see Definition \ref{defi: canonical form}. Then 
 $\phi_{T^*_X}(\theta_{X'})$ is the class of the algebra of differential operators ${F_X}_*D_X$ in the Brauer group $Br(T^*_{X'}).$
\end{prop}

\begin{proof} This follows from \cite[Propositions 4.4 and 4.2]{Ogus-Vologodsky}. (See also \cite[Proposition 3.11]{BezrBr}.)
\end{proof}

\section{Lagrangianity} \label{sec:lagrangianity}

In this section, we complete the proof of Theorem \ref{thm:main thm}.

\subsection{Nice compactification of the $p$-supports} \label{subsec:compactification}
Here for $p$ large enough, we show the existence of normal-crossings compactifications of the regular locus of the $p$-supports which are well-behaved as $p$ varies. The poles of the canonical form at infinity of this nice compactification have nice properties. This will allow us to transfer characteristic zero results to $p$-supports.

\subsubsection{Poles}

To fix notations, we first recall some definitions.  
Let $S$ be a scheme, let $\overline{Y}$ be a smooth $S$-scheme and let $D\subset \overline{Y}$ be a divisor with normal crossings relative to $S,$ defined by the invertible ideal sheaf $I.$ 
Denote by $Y\xhookrightarrow{j}\overline{Y}$ the affine open embedding complement to $D.$  
For an $\mathcal{O}_{\overline{Y}}$-module $\mathcal{F}$ and $n\in \mathbb{Z},$ we let $\mathcal{F}(nD)$ denote $\mathcal{F}\otimes_{\mathcal{O}_{\overline{Y}}}\mathcal{I}^{\otimes_{\mathcal{O}_{\overline{Y}}}(-n)}.$ We have $j_*\Omega^m_{Y/S}= \bigcup_{n\geq0}\Omega^m_{\overline{Y}/S}(nD),$ where $\Omega^m_{Y/S}$ is the sheaf of relative differential forms of degree $m$ on $Y/S.$ Note that if $S$ is the spectrum of a field and when there is no risk of confusion, we will (and have already) commit the abuse of notation to denote $\Omega^m_{Y/S}$ by $\Omega^m_{Y}.$

\begin{defi} A local section of $j_*{\Omega}^m_{Y/S}$ which is in ${\Omega}^m_{\overline{Y}/S}(nD)$ is said to have poles of order at most $n$ along $D.$
\end{defi}

We denote the logarithmic de Rham complex by $(\Omega^{\bullet}_{\overline{Y}/S}(logD), d),$ see \cite[II \S 3]{DeligneLNM}. It is a subcomplex of $j_*(\Omega^{\bullet}_{Y/S}, d_{Y/S}).$ 

\begin{defi} A local section of $j_*\Omega^m_{Y/S}$ which is in $\Omega^m_{\overline{Y}/S}(logD)$ is said to have logarithmic poles along $D.$
\end{defi}

\subsubsection{Hilbert scheme}

Let us fix coordinates $\{x_1, \dots, x_n\}$ on $\mathbb{A}^n_{\mathbb{Z}}.$  
They induce by base-change coordinates $\{x_1, \dots, x_n\}$ on $\mathbb{A}^n_S$ and $\{x_1, \dots, x_n, y_1, \dots, y_n\}$ on its relative cotangent space $T^*_{\mathbb{A}^n_S},$ for any scheme $S.$ We thus have a filtration of $\mathcal{O}_{T^*_{\mathbb{A}^n_S}}$ by the order in $\{x_1, \dots, x_n, y_1, \dots, y_n\}$ and an open embedding $T^*_{\mathbb{A}^n_S} \hookrightarrow \mathbb{P}^{2n}_S,$ by the Rees construction (Proposition \ref{prop: Rees}). 
If $S$ is the spectrum of a field $K$ of positive characteristic, we also deduce identifications of the base-change by Frobenius $(\mathbb{A}^n_K)' = \mathbb{A}^n_K$ and $T^*_{(\mathbb{A}^n_K)'}= T^*_{\mathbb{A}^n_K}.$ This embedding is the same as that considered on page \pageref{page:Rees}. 

Let $P\in \mathbb{Q}[X]$ be a rational polynomial and let $\mathcal{H}:= \mathcal{H}ilb^P_{\mathbb{P}^{2n}_{\mathbb{Z}}}$ be the Hilbert scheme of $\mathbb{P}^{2n}_{\mathbb{Z}}$ of index $P,$ see \cite[p.17]{FGA Hilb}. We let $\mathcal{Z}_{\mathcal{H}} \subset\mathbb{P}^{2n}_{\mathbb{Z}}\times_{\Spec(\mathbb{Z})}\mathcal{H}$ be the associated universal flat closed subscheme. The following definition is convenient:

\begin{defi} \label{defi: nice}
A (locally closed) subscheme $\mathcal{S} \xhookrightarrow{i} \mathcal{H}$ is nice if
\begin{enumerate}
\item $\mathcal{S}$ is irreducible and generically of characteristic $0,$

\item Let $\mathcal{Z}$ be the restriction $\mathcal{Z}_{\mathcal{H}}\times_{\mathcal{H}}\mathcal{S}$ of $\mathcal{Z}_{\mathcal{H}}$ to $\mathcal{S}.$ There exists an open $U\subset \mathcal{Z}^{red}$ surjecting onto $\mathcal{S}$ and an open immersion $U\cap T^*_{\mathbb{A}^n_{\mathcal{S}}}=: \mathcal{Y} \xhookrightarrow{j} \overline{\mathcal{Y}}$ such that $\overline{\mathcal{Y}}$ is a smooth projective $\mathcal{S}$-scheme and the complement $\Delta$ of $j$ is a divisor with normal crossings relative to $\mathcal{S},$

\item \label{item: vanishing}
Let $\mu$ be the generic point of $\mathcal{S}.$ Then $(d\theta_{\mathbb{A}^n_{\mathcal{S}}}|_\mathcal{Y})= 0$ if and only if $(d\theta_{\mathbb{A}^n_{\mathcal{S}}}|_\mathcal{Y})_\mu=0,$

\item \label{item: log pole}
The restriction $d\theta_{\mathbb{A}^n_{\mathcal{S}}}|_\mathcal{Y}$ of the symplectic form $d\theta_{\mathbb{A}^n_{\mathcal{S}}}$ to $\mathcal{Y}$ has logarithmic poles along $\Delta$ if and only if there exists a closed point $s\in \mathcal{S}$ such that the fiber $(d\theta_{\mathbb{A}^n_{\mathcal{S}}}|_\mathcal{Y})_s$ has logarithmic poles along $\Delta_s.$

\end{enumerate}
\end{defi}

We may now state and prove the existence of compactification result. 

\begin{prop}\label{prop:compactification}
There exists a positive integer $N$ such that $\mathcal{H}[\frac{1}{N}]$ has a finite partition into nice subschemes.
\end{prop}

\begin{proof} Let $X\hookrightarrow \mathcal{H}$ be a subscheme. Denote $P(X)$ the property that there exists an integer $N_X>0$ such that $X[\frac{1}{N_X}]$ is empty or admits a partition into nice subschemes of $\mathcal{H}.$ We have $P(X)$ if and only if $P(X^{red})$ holds. Let us prove $P(X)$ holds for all reduced closed subschemes $X$ of $\mathcal{H},$ by Noetherian induction.

Let $S\hookrightarrow \mathcal{H}$ be a reduced subscheme. By Noetherian induction, we assume $P(T)$ holds for all reduced closed subschemes $T$ of $S,$ distinct from $S.$ Clearly there is an integer $N>0$ such that the irreducible components of $S[\frac{1}{N}]$ have all generic points of characteristic zero. If $S[\frac{1}{N}]$ is empty then we have $P(S).$ We thus from now on assume that $S[\frac{1}{N}]$ is not empty.

Let $z\in S[1/N]$ be the generic point of an irreducible component and let $\mathcal{Z}_z$ be the fiber of $(\mathcal{Z}_{\mathcal{H}}\times_{\mathcal{H}}S[1/N])^{red}$ at $z.$ Then $\mathcal{Z}_z\subset\mathbb{P}^{2n}_{k(z)}$ is a scheme of finite type over the field of characteristic zero $k(z).$ It is reduced by \cite[Proposition 8.7.2 a)]{EGAIV} and hence contains an open dense smooth subset $U\subset\mathcal{Z}_z.$  
By the resolution of singularities in characteristic zero, there is an open immersion $Y \overset{j}\hookrightarrow \overline{Y}$ of the quasi-projective variety $Y:=U\cap T^*_{\mathbb{A}^n_{k(z)}}\subset\mathbb{P}^{2n}_{k(z)}$ into a smooth projective scheme $\overline{Y}$ over $k(z)$ which is the complement of a divisor $D$ with normal crossings relative to $k(z).$
Hence by \cite[Th\'eor\`eme 8.10.5 and Proposition 17.7.8]{EGAIV}, there are an open affine neighborhood $\mathcal{T}$ of $z,$ a non-empty smooth open subset $\mathcal{U}\subset(\mathcal{Z}_{\mathcal{H}}\times_{\mathcal{H}}\mathcal{T})^{red}$ and an open $\mathcal{T}$-immersion $\mathcal{U}\cap T^*_{\mathbb{A}^n_\mathcal{T}}=:\mathcal{Y} \overset{j}\hookrightarrow \overline{\mathcal{Y}}$ into a smooth projective $\mathcal{T}$-scheme which is the complement of a divisor $\mathcal{D}$ with normal crossings relative to $\mathcal{T}.$ Moreover, since smooth morphisms are open, $\mathcal{T}$ may be chosen such that $\mathcal{U}$ surjects on $\mathcal{T}.$ Finally by \cite[Proposition 2.1.9 (ii)]{EGAI} it can also be chosen to be integral and hence irreducible, since $S$ is reduced.

Let $\theta$ be the restriction to $\mathcal{Y}$ of the canonical form $\theta_{\mathbb{A}^n_{\mathcal{T}}/\mathcal{T}}$ on $T^*_{\mathbb{A}^n_{\mathcal{T}}}.$ If $d\theta$ vanishes on the generic fiber $\mathcal{Y}_z,$ then there is a dense open subset $V\subset \mathcal{T}$ above which $d\theta$ vanishes.  
If $d\theta$ does not vanish on the generic fiber of $\mathcal{Y},$ then we set $V =\mathcal{T}.$

Finally, there is a dense open subset $W\subset V$ such that $d\theta|_{\mathcal{Y}|_W}$ has logarithmic poles along $\mathcal{D}|_W$ as soon as there is a fiber on which it has logarithmic poles. Indeed there is a closed subset $\mathcal{L}$ of $V$ such that for all closed points $v \in V,$ $d\theta|_{\mathcal{Y}_v}$ has logarithmic poles along $\mathcal{D}_v$ if and only if $v\in \mathcal{L}.$ We prove it as follows. There is an integer $m>0$ such that $d\theta|_{\mathcal{Y}}$ is a section of the locally free $\mathcal{O}_{\overline{\mathcal{Y}}}$-module ${\Omega}^2_{\overline{\mathcal{Y}}/V}(m\mathcal{D}).$ Let us locally choose a basis $\{b_1, \dots, b_j\}$ of ${\Omega}^2_{\overline{\mathcal{Y}}/V}(m\mathcal{D})$ extending a basis $\{b_1, \dots, b_l\}$ of ${\Omega}^2_{\overline{\mathcal{Y}}/V}(log\mathcal{D}) \subset {\Omega}^2_{\overline{\mathcal{Y}}/V}(m\mathcal{D})$ such that  for each closed point $v\in V,$ the restriction $\{b_1|_{\mathcal{Y}_v}, \dots, b_j|_{\mathcal{Y}_v}\}$ is a basis of ${\Omega}^2_{\overline{\mathcal{Y}}_v}(m\mathcal{D}_v)$ extending a basis $\{b_1|_{\mathcal{Y}_v}, \dots, b_l|_{\mathcal{Y}_v}\}$ of ${\Omega}^2_{\overline{\mathcal{Y}}_v}(log\mathcal{D}_v).$ That such a basis exists follows from the case of $(\mathbb{A}^d_V, \{x_1\dots x_r = 0\}),$ which is immediate. Locally we thus have $d\theta|_{\mathcal{Y}}= \Sigma_{i=1}^{i=j} \alpha_ib_i$ for some unique sections $\alpha_1, \dots, \alpha_j$ of $\mathcal{O}_{\overline{\mathcal{Y}}}.$ Moreover its restriction $d\theta|_{\mathcal{Y}_v}$ has logarithmic poles along $\mathcal{D}_v$ if and only if the restrictions $\alpha_{l+1}|_{\overline{\mathcal{Y}}_v}, \dots, \alpha_j|_{\overline{\mathcal{Y}}_v}$ vanish. Let $O \subset \overline{\mathcal{Y}}$ be the open complement of the closed subset of $\overline{\mathcal{Y}}$ locally given by $\{\alpha_{l+1}=\dots=\alpha_j=0\}.$ For each closed point $v\in V,$ we thus have that $d\theta|_{\mathcal{Y}_v}$ has logarithmic poles along $\mathcal{D}_v$ if and only if $O\cap \overline{\mathcal{Y}}_v$ is empty, that is $v\not\in\pi(O),$ where $\pi$ is the structure map $\overline{\mathcal{Y}}\to V.$ But $\pi(O)$ is open since $\pi$ is smooth and thus open. We have that $\mathcal{L}$ is the complement of $\pi(O).$
If $\mathcal{L}\subsetneq V,$ then we let $W$ be the open complement $V\setminus\mathcal{L}.$ Otherwise we set $W=V.$

Thus we have found a nice non-empty open subset $W\subseteq S.$ Since by Noetherian induction we have $P(S\setminus W),$ it follows immediately that we have $P(S).$ This concludes the proof of the proposition. 
\end{proof}

\subsection{Action of the $p$-curvature operator on the order of poles} \label{subsec: poles}
Let $\eta$ be a 1-form in the image of the $p$-curvature operator. Here we will show that if $\eta$ has poles of small order (i.e. $\leq p-1$), then its exterior derivative $d\eta$ has at most simple poles and hence has logarithmic poles. The following definition is convenient.

\begin{defi} \label{defi: good coordinates}
Let $D$ be a divisor in a smooth variety $X$ over a field $K.$ Etale coordinates $\{x_1, \dots, x_n\}: V_x \to \mathbb{A}^n_K$ in the neighbourhood $V_x$ of a point $x\in X$ are $D$-good, or good if there is no risk of confusion, if there is $0\leq r\leq n$ such that $D\cap V_x=\{x_1\dots x_r =0\}.$
\end{defi}
Thus $D\subset X$ has normal crossings if and only if $X$ has an \'etale covering $\pi$ admitting $\pi^{-1}(D)$-good \'etale coordinates in the neighbourhood of each of its points. 

Throughout this subsection, $\overline{Y}$ will denote a smooth variety over a perfect field $K$ of positive characteristic $p$ and $D\subset \overline{Y}$ a normal crossings divisor of open complement $Y\xhookrightarrow{j} \overline{Y}.$ 
As above, we use ' to denote the base-change by the Frobenius automorphism of $K.$ Recall Definition \ref{defi: p-curvature} of the $p$-curvature operator $F_*Z^1(\Omega^{\bullet}_Y) \xrightarrow{W^{\star}-C_Y} \Omega^{1}_{Y'}.$
Here is the main result of the subsection:
\begin{prop} \label{prop:log poles} Let $\mathcal{I}m(W^{\star}-C_Y)$ be the image of abelian sheaves for the Zariski topology and let $\mathcal{I}$ be the intersection of the subsheaves ${\Omega}^1_{\overline{Y}'}((p-1)D')$ and $j'_*\mathcal{I}m(W^{\star}-C_Y)$ of $j'_*\Omega^{2}_{Y'}.$ Then $d(\mathcal{I}) \subset \Omega^2_{\overline{Y}'}(logD').$
\end{prop}

\begin{proof} The assertion is \'etale local on $\overline{Y}.$ Hence by definition of divisor with normal crossings, we may and do assume that there are good \'etale coordinates in the neighbourhood of every point of $\overline{Y}.$ In this situation, let us prove that $\mathcal{I}\subseteq j'_*B^1\Omega^\bullet_{Y'} + \Omega^1_{\overline{Y}'}(log D'),$ with $B^1\Omega^\bullet_{Y'}$ the exact 1-forms. This immediately implies the proposition.

Let $\eta$ be a local section of $\mathcal{I}.$ By Lemma \ref{lmm:j and Im commute}, we may assume that locally, there is a local section $\rho\in j'_*F_*Z^1(\Omega^{\bullet}_Y)$ such that $\eta= j'_*(W^{\star}-C_Y)\rho.$ Moreover, by $p^{-1}$-linearity of the Cartier operator $C_Y,$ see  \cite[Lemme 0.2.5.4]{Illusie dR-W}, $\rho$ has poles of order at most $p-1$ along $D.$

Hence by Lemma \ref{lmm: poles}, $\rho$ is locally a sum of an exact form and a closed form with logarithmic poles along $D.$ The proposition follows since, by \cite[(7.2.4) of Theorem 7.2]{Katz}, the Cartier operator $C_Y$ preserves forms with logarithmic poles.
\end{proof}

\begin{lmm} \label{lmm: poles} 
Suppose that there are $D$-good \'etale coordinates in the neighbourhood of each point of $\overline{Y}.$
Then there is an inclusion of quasi-coherent subsheaves of $j'_*F_*\Omega^1_Y:$
$$j'_*Z^1(F_*\Omega^\bullet_Y) \cap F_*\Omega^1_{\overline{Y}}((p-1)D) \subseteq j'_*B^1(F_*\Omega^\bullet_Y) + Z^1(F_*\Omega^\bullet_{\overline{Y}}(log D)),$$ where $B^1(F_*\Omega^\bullet_Y):= \mathrm{Im}(F_*\mathcal{O}_Y \xrightarrow{F_*d} F_*\Omega^1_Y)$ and $Z^1(F_*\Omega^\bullet_{\overline{Y}}(log D)):= \ker(F_*\Omega^1_{\overline{Y}}(log D)\xrightarrow{F_*d}F_*\Omega^2_{\overline{Y}}(log D)).$  
\end{lmm}

\begin{proof} 
The assertion being local, we may assume that we have global good \'etale coordinates $\overline{Y}\xrightarrow{\{y_1,\dots, y_n\}} \mathbb{A}^n_K.$ Moreover it is clear that the lemma reduces to $\overline{Y}= \mathbb{A}^n_K$ by pullback. So we may assume that $(\overline{Y}, D)= (\mathbb{A}^n_K, \{y_1\dots y_r =0\})$ for some $0\leq r\leq n,$ where $\{y_1,\dots, y_n\}$
are coordinates on $\mathbb{A}^n_K.$ 

The coordinates induce a splitting of the canonical short exact sequence
$0\to B^1(F_*\Omega^{\bullet}_Y)\to Z^1(F_*\Omega^{\bullet}_Y)\xrightarrow{C_Y}{\Omega}^1_{Y'} \to 0$ by ${\Omega}^1_{Y'}\xrightarrow{\delta} Z^1(F_*\Omega^{\bullet}_Y),$ sending a local section $\sum\limits_{i=1}^{n} a_idy_i$ to $\delta(\sum\limits_{i=1}^{n} a_idy_i)= \sum\limits_{i=1}^{n} a_i^py_i^{p-1}dy_i,$ where we have identified $\mathbb{A}^n_K$ and its base-change by Frobenius using the coordinates. Applying $j'_*$ we get that a local section $\rho$ of $j'_*Z^1(F_*\Omega^\bullet_Y)$ may uniquely be written as a sum $\rho= \rho_b+\rho_c,$ where $\rho_b$ is a local section of $j'_*B^1(F_*\Omega^\bullet_Y),$ and 
$\rho_c= \sum\limits_{i=1}^{n} a_i^py_i^{p-1}dy_i,$ with $\{a_1, \dots, a_n\}$ uniquely determined local sections of $j'_*\mathcal{O}_{Y'}.$

Let us show that if $\rho$ is a local section of $j'_*Z^1(F_*\Omega^\bullet_Y) \cap F_*\Omega^1_{\overline{Y}}((p-1)D),$ then $\rho_c$ is a local section of $Z^1(F_*\Omega^\bullet_{\overline{Y}}(log D)).$ This implies the lemma.
First one sees immediately by a computation in coordinates that if $\rho_c$ has poles of order at most $p-1,$ then it actually has poles of order at most $1.$ Thus it has logarithmic poles along $D,$ since it is closed.

We conclude the proof by showing that if $\rho$ has poles of order at most $p-1,$ then so does $\rho_c$ (and hence $\rho_b$).
Set $h:= y_1\dots y_r.$ We claim that $(h^n\rho)_c= h^n\rho_c,$ for all $n\geq0,$ where the multiplication by $h^n$ denotes here the action on  
$\Omega^1_Y.$ Since for all local sections $\xi$ of $j'_*Z^1(F_*\Omega^\bullet_Y), \xi_c= j'_*\delta\circ j'_*C_Y(\xi)$ and $j'_*\delta\circ j'_*C_Y$ preserves regular forms (i.e. forms having poles of order at most $0$ along $D$), the claim indeed implies the assertion. But the claim is easily checked by adjoining to the domain $\mathcal{O}(Y)$ all the $p$-th roots of its elements. Indeed the composition of the $p^{-1}$-linear morphism $j'_*C_Y$ and the $p$-linear morphism $j'_*\delta$ is linear. 
\end{proof}

\begin{lmm} \label{lmm:j and Im commute}

Let $\overline{X}$ be a smooth variety over a perfect field $K$ of positive characteristic $p$ and let $D\subset \overline{X}$ be a normal crossings divisor of open complement $X\xhookrightarrow{j} \overline{X}.$
Then the canonical inclusion $\mathcal{I}m(j'_*(W^{\star}-C_X)) \hookrightarrow j'_*\mathcal{I}m(W^{\star}-C_X)$ is an isomorphism.

\end{lmm}
\begin{proof} The $p$-curvature exact sequence of abelian sheaves on $X'$ for the Zariski topology $0\to \mathcal{O}_{X'}^\times \xrightarrow{F^*} F_*\mathcal{O}_X^\times \xrightarrow{F_*dlog} Z^1(F_*\Omega^{\bullet}_X) \xrightarrow{W^{\star}-C_X} \Omega^{1}_{X'},$ see  \cite[Corollaire 0.2.1.18]{Illusie dR-W}, breaks into two short exact sequences: $0\to \coker F^* \xrightarrow{F_*dlog} Z^1(F_*\Omega^{\bullet}_X) \xrightarrow{W^{\star}-C_X} \mathcal{I}m(W^{\star}-C_X) \to 0$ and $0\to \mathcal{O}_{X'}^\times \xrightarrow{F^*} F_*\mathcal{O}_X^\times \to \coker F^* \to 0.$ The long exact sequence for $Rj'_*$ of the former:
\begin{multline*}
0\to j'_*\coker F^* \xrightarrow{j'_*F_*dlog} j'_*Z^1(F_*\Omega^{\bullet}_X) \xrightarrow{j'_*(W^{\star}-C_X)} j'_*\mathcal{I}m(W^{\star}-C_X) \to R^1j'_*\coker F^*\to...
\end{multline*}
implies that the lemma follows from the vanishing of $R^1j'_*\coker F^*.$ 

This in turn would follow from the vanishing of $R^1j'_*(F_*\mathcal{O}_X^\times)$ and $R^2j'_*\mathcal{O}_{X'}^\times,$ by the long exact sequence for $Rj'_*$ of the other short exact sequence:
 
$$...\to R^1j'_*\mathcal{O}_{X'}^\times \to R^1j'_*(F_*\mathcal{O}_X^\times) \to R^1j'_*\coker F^* \to R^2j'_*\mathcal{O}_{X'}^\times\to...$$  

But the direct image $F_*$ preserves flasque sheaves and is exact, since $F$ is a homeomorphism. Hence by \cite[III Corollary 8.3]{HartshorneAG}, $R^qF_*(G)=0$ for all abelian sheaves $G$ and all $q>0.$ Thus $R^1j'_*(F_*\mathcal{O}_X^\times) \cong R^1(j'_*\circ F_*)(\mathcal{O}_X^\times)=R^1(F_* \circ j_*)(\mathcal{O}_X^\times) \cong F_*R^1j_*(\mathcal{O}_X^\times).$ 
We have therefore shown that the Lemma is implied by the vanishing of $R^1j_*(\mathcal{O}_X^\times)$ and $R^2j'_*\mathcal{O}_{X'}^\times.$ This follows from Lemma \ref{lmm:vanishing of higher direct images}
below, since smooth schemes are locally factorial.\end{proof}

\begin{lmm} \label{lmm:vanishing of higher direct images}

Let $U \overset{j}\hookrightarrow Y$ be an open immersion. Suppose that $Y$ is a locally factorial Noetherian scheme. Then $R^qj_*(\mathcal{O}_U^\times)=0,$ for all $q>0.$

\end{lmm}
\begin{proof} By \cite[Proposition 8.1 of Chapter III]{HartshorneAG}, $R^qj_*(\mathcal{O}_U^\times)$ is the abelian sheaf associated to the presheaf $V \mapsto H^q(U\cap V,\mathcal{O}_{U\cap V}^\times),$ for all $V$ open in $Y.$ But by Lemma \ref{lmm:vanishing of Zariski cohomology}, $H^q(U\cap V,\mathcal{O}_{U\cap V}^\times)=0$ for all $q\geq 2.$ Thus $R^qj_*(\mathcal{O}_U^\times)=0,$ for all $q\geq 2.$  

For $q=1,$ $R^1j_*(\mathcal{O}_U^\times)$ is the abelian sheaf associated to the presheaf of Picard groups: $V \mapsto 
Pic(U\cap V),$ for all $V$ open in $Y.$ Let $\mathcal{L}_{U\cap V} \in Pic(U\cap V).$ By \cite[Corollaire 21.6.11]{EGAIV}, it extends to an invertible sheaf $\mathcal{L}_V \in Pic(V).$ Hence $\mathcal{L}_{U\cap V}$ is trivial on the restriction to $U$ of an open covering of $V$ trivializing $\mathcal{L}_V.$ Thus the corresponding section of the associated sheaf $R^1j_*(\mathcal{O}_U^\times)$ vanishes locally, hence is $0.$ This implies that $R^1j_*(\mathcal{O}_U^\times)$ vanishes. 
\end{proof}

\subsection{Conclusion} \label{subsec:conclusion}

We now combine most of the results above to prove our main theorem. 
\begin{proof}[Proof of Theorem \ref{thm:main thm}:] 
 
By Theorem \ref{thm:purity}, there is a dense open subset $U_1$ of $S$ such that $p$-supp$(M_u)$ is equidimensional of dimension $n,$ for all closed points $u$ of $U_1.$

By Proposition \ref{prop:reduction to affine space}, we may assume that $X/S=\mathbb{A}^n_S/S.$ Thus by Theorem \ref{thm:splitting in affine space case}, there is a dense open subset $U_2$ of $U_1$ such that the Azumaya algebra $F_*D_{\mathbb{A}^n_{k(u)}}$ splits on the regular locus of each irreducible component of the $p$-support of $M_u,$ for every closed point $u$ of $U_2.$ Let $\theta_u= \theta_{(\mathbb{A}^n_{k(u)})'/\Spec k(u)}$ be the canonical form on $T^*_{(\mathbb{A}^n_{k(u)})'}.$ Then the restriction of $\theta_u$ to the regular locus $Z_u^{reg}$ of an irreducible component $Z_u$ of $p$-supp$(M_u)$ is in the image of the $p$-curvature operator $\mathcal{I}m(W^{\star}-C_{Z_u^{reg}}),$ where we have identified $Z_u^{reg}$ with ${Z_u^{reg}}',$ as $k(u)$ is perfect. This follows from the description of $\ker\overline{\phi_Y}$ in Proposition \ref{prop:kernel of phi}, since $\phi_{T^*_{\mathbb{A}^{n}_{k(u)}}}(\theta_u)$ is the class of $F_*D_{\mathbb{A}^n_{k(u)}}$ in the Brauer group by Proposition \ref{prop:canonical form maps to algebra of differential operators} and $\phi_Y$ is commutes with the pullback by Lemma \ref{lmm:functoriality of phi}.

By Theorem \ref{thm:bound} and using its notations, there are a dense open subset $U_3$ of $U_2$ and an integer $e>0$ such that for all closed points $u$ of $U_3,$ the degree of each irreducible component $Z_u$ of the $p$-support of $M_u$ is not larger than $e.$ Hence the Hilbert polynomial of $Z_u$ belongs to a finite set $\{P_1, \dots, P_r\},$ independent of $u.$ Indeed the Hilbert polynomial of $Z_u$ is the same as that of $Z_u \times_{\Spec k(u)}\Spec\overline{k(u)}$ for an algebraic closure $\overline{k(u)}$ of the residue field $k(u).$ We can then conclude using Chow coordinates, see \cite[Lemme 2.4 and Th\'eor\`{e}me 2.1(b)]{FGA Hilb}. Moreover by Proposition \ref{prop:compactification}, there is an integer $N>0$ such that the localization $\coprod_{i=1}^{i=r}\mathcal{H}ilb^{P_i}_{\mathbb{P}^{2n}_{\mathbb{Z}}}[\frac{1}{N}]$ of the Hilbert scheme of index  
$\{P_1, \dots, P_r\}$ has a finite partition into nice subschemes, see Definition \ref{defi: nice}. Hence, in particular, there are a dense open subset $U_4$ of $U_3$ and an integer $N'>0$ satisfying the following. For all closed points $u$ of $U_4,$ each irreducible component $Z_u$ of the $p$-support of $M_u$ has a smooth dense open subset $Y_u\subset Z_u$ embedding as an open into a smooth projective variety $Y_u \xhookrightarrow{j_u}\overline{Y_u}$ with complement a divisor $D_u$ with normal crossings. Moreover the restriction of the canonical form $\theta_u|_{Y_u}$ has poles of order at most $N'$ along $D_u.$ 

Hence by further inverting a positive integer, there exists an open dense subset $U_5$ of $U_4$ such that for all closed points $u$ of $U_5, \theta_u|_{Y_u}$ has poles of order at most $p_u-1$ along $D_u,$ where $p_u$ is the characteristic of $k(u).$ Thus since $\theta_u|_{Y_u}$ is in the image of the $p$-curvature operator, Proposition \ref{prop:log poles} implies that its exterior derivative $d\theta_u|_{Y_u}$ has logarithmic poles along $D_u.$

By construction of $\overline{Y_u}$ there are an irreducible $U_5$-scheme $H$ generically of characteristic zero, a smooth projective $H$-scheme $\overline{Y}$ and a divisor $\Delta$ of $\overline{Y}$ with normal crossings relative to $H$ such that the complement $Y:= \overline{Y}-\Delta$ is identified with an open in $T^*_{\mathbb{A}^n_H}.$ And if we denote by $Y \xhookrightarrow{j} \overline{Y}$ the open embedding, there is a closed point $t$ of $H$ such that $j_u$ (resp. $\theta_u|_{Y_u}$) is the specialization of $j$ (resp. $\theta_{\mathbb{A}^n_H/H}|_Y$) at $t.$ By point (\ref{item: log pole}) of the definition of nice subscheme, the above implies that $d\theta_{\mathbb{A}^n_H/H}|_Y$ has logarithmic poles along $\Delta.$ Hence so does its restriction $d\theta_{\mathbb{A}^n_{k(\gamma)}}|_{Y_{\gamma}}$ to the fiber at the generic point $\gamma$ of $H.$ Clearly, being globally exact, the class of the symplectic form $d\theta_{\mathbb{A}^n_{k(\gamma)}}|_{Y_{\gamma}}$ vanishes in the hypercohomology of the de Rham complex $({j_{\gamma}}_*\Omega_{Y_{\gamma}}^\bullet, {j_{\gamma}}_*d).$ But since the residue field of $\gamma$ is of characteristic zero, the canonical inclusion of complexes $\Omega_{\overline{Y_{\gamma}}}(log \Delta_{\gamma}) \subset {j_{\gamma}}_*\Omega_{Y_{\gamma}}^\bullet$ is a quasi-isomorphism by \cite[Proposition 3.1.8]{Hodge II}. Hence the class of $d\theta_{\mathbb{A}^n_{k(\gamma)}}|_{Y_{\gamma}}$ vanishes in the hypercohomology of the logarithmic de Rham complex, as well. Thus by the degeneracy at $E_1$ of the logarithmic Hodge to de Rham spectral sequence (\cite[Corollaire 3.2.13 (ii) and Corollaire 3.2.14]{Hodge II}), $d\theta_{\mathbb{A}^n_{k(\gamma)}}|_{Y_{\gamma}}=0.$ Therefore by point (\ref{item: vanishing}) of the definition of nice subscheme, $d\theta_{\mathbb{A}^n_H/H}|_Y=0.$ It follows directly that for all closed points $u$ of $U_5, d\theta_u|_{Y_u}=0.$ This concludes the proof of the theorem.
\end{proof}
Let us  
explain how our results 
imply the involutivity of the singular support of a holonomic $\mathcal{D}$-module $M.$ The latter is an important special case of the involutivity of the singular support of an arbitrary nonzero coherent $\mathcal{D}$-module, which was originally proved in \cite{KKS}. A purely algebraic proof was later given in \cite{Gabber}. 
\begin{cor} \label{cor: Gabber} Let Y be a smooth variety over a field $L$ of characteristic zero and let $M$ be a nonzero holonomic left $\mathcal{D}_Y$-module. Then the singular support of $M$ is a Lagrangian subvariety of $(T^*_Y, \omega_Y).$

\end{cor}

\begin{proof} The corollary immediately reduces to an equidimensional and affine $Y.$ One proves that the singular support is equidimensional of dimension $\dim Y$ by the methods of Section \ref{sec:dimension of the p-supports}, purely in characteristic zero, see e.g. \cite[A:IV Theorem 5.2]{Bjork2}. 

Let us now show that the symplectic form $\omega_Y$ vanishes on the regular locus of the singular support $SS(M).$ Consider a spreading out  
$(X, N)$ of $(Y, M).$ Namely, there are a finitely generated subring $A\subset L,$ a smooth affine scheme $X$ over $S:=\Spec(A)$ and a coherent left $D_{X/S}$-module $N$ such that, for $\mu$ the generic point of $S,$ there are compatible isomorphisms $L\otimes_{k(\mu)}X_\mu\cong Y$ and $L\otimes_{k(\mu)}N_\mu\cong M.$ Let us define the closed subset $SS(N)\subset T^*_{X/S}$ to be the support of the associated graded $gr^\Gamma N,$ where $\Gamma$ is a good filtration of $N,$ see Definition \ref{defi: good filtration}. As usual the subset $SS(N)$ does not depend on the choice of good filtration of $N.$ We claim that there is a dense open subset $S^*$ of $S$ such that, for all closed points $u$ of $S^*,$ the symplectic form $\omega_{X_u}$ vanishes on the regular locus of the fiber $SS(N)_u$ of $SS(N).$ Since $SS(M)= L\otimes_{k(\mu)}{SS(N)}_\mu,$ the corollary follows from the claim.

By generic freeness of $gr^\Gamma N$ on $S,$ there is a dense open subset $S_1$ of $S$ such that the fiber $\Gamma_u$ of $\Gamma$ is a good filtration and $SS(N)_u=SS(N_u),$ for all closed points $u$ of $S_1.$ Moreover it is clear that Lemma \ref{lmm:pGamma is good} and its proof generalize to the Bernstein filtration replaced by the filtration by the order of differential operators, and $\mathbb{A}^n_{k(u)}$ replaced by $X_u.$ Thus (\ref{eq: grZ}) gives an isomorphism of $gr^\Phi Z(D_{X_u})$-modules $gr(gr^{\Phi(\Gamma_u)} r_*N_u) \simeq F_*gr^{\Gamma_u} N_u,$ where $\Phi$ is the filtration induced by the order of differential operators on the center $Z(D_{X_u})$ of $D_{X_u}$ and we have used the notations of the proof of Lemma \ref{lmm:pGamma is good}. In addition, we have that the symplectic form $\omega_{X_u}$ vanishes on the regular locus of the support of $gr^{\Gamma_u} N_u$ if and only if $\omega_{X'_u}$ vanishes on the regular locus of the support of $F_*gr^{\Gamma_u} N_u.$ Indeed these regular loci and symplectic forms are mapped to one another by the projection map $T^*_{X'_u}\to T^*_{X_u},$ which is an isomorphism since $k(u)$ is perfect.
Hence the symplectic form $\omega_{X_u}$ vanishes on the regular locus of the support of $gr^{\Gamma_u} N_u$ if and only if $\omega_{X'_u}$ vanishes on the regular locus $T_u^{reg}$ of the support $T_u$ of $gr^{\Phi(\Gamma_u)} r_*N_u.$ 

But using the Rees module with respect to $\Phi,$ we see that $T_u$ is the reduced scheme associated to the fiber at the origin of a flat $\mathbb{A}^1$-scheme of generic fiber isomorphic to the $p$-support of $N_u.$ And by Theorem \ref{thm:main thm}, there is a dense open subset $S_2$ of $S_1$ such that the $p$-support of $N_u$ is Lagrangian, for all closed points $u$ of $S_2.$
Thus the claim would follow if the reduced scheme associated to the fiber at the origin of a flat $\mathbb{A}^1$-scheme of Lagrangian generic fiber was Lagrangian. This is not in general true in positive characteristic.  
In our case however,  
the proof reduces to the case $X_u = \mathbb{A}^n$ and we may use Theorem \ref{thm:bound}. Hence there is a dense open subset $S_3$ of $S_2$ such that for all closed $u\in S_3,$ the generic fiber of the $\mathbb{A}^1$-scheme is Lagrangian and of degree bounded independently of $u\in S_3$ (for the projective embedding of Theorem \ref{thm:bound}). We may thus use the Hilbert scheme to conclude that there is an integer $l>0$ such that for all $u\in S^*:= S_3[\frac{1}{l}],$  
the assertion that the reduced scheme associated to the fiber at the origin of the above $\mathbb{A}^1$-scheme is Lagrangian, follows from the case of characteristic zero. But the latter is a special case of Lemma \ref{lmm: symplectic specialization} below, applied to the symplectic form.  
This concludes the proof of the claim and thus of the corollary.\end{proof}

\begin{lmm} \label{lmm: symplectic specialization}
	Let $S$ be an integral scheme over a field $L$ of characteristic $0.$ Let $\frak{M}$ be a smooth $S$-scheme, let $Z \subset \frak{M}$ be a subscheme which is a reduced, flat and surjective $S$-scheme of relative dimension $d,$ and let $\alpha$ be a relative differential form on $\frak{M}$ over $S.$ Assume that the restriction of $\alpha$ to the smooth locus of the generic fibre of $Z$ vanishes and that $S$ is the spectrum of a Dedekind ring. Then the restriction of $\alpha$ to the smooth locus of the reduced scheme associated to every closed fibre of $Z$ vanishes.\end{lmm}

\begin{proof}
We will reduce the proof to the case of a smooth $Z,$ which is immediate. Indeed if $Z$ is smooth over $S,$ the points $s$ of the base $S$ such that the restriction of $\alpha$ vanishes on the fibre $Z_s$ form a closed subset $V_\alpha \subset S.$ Since the generic point of $S$ belongs to $V_\alpha$ by hypothesis, we have $V_\alpha=S,$ which is what we wanted to prove.

By the Reduced Fibre Theorem \cite[\href{https://stacks.math.columbia.edu/tag/09IL}{Theorem 09IL}]{Stacks}, there is an integral affine scheme $T$ and a finite surjective morphism $T \xrightarrow{\phi} S,$ such that if we denote by $Z_T$ the base-change of $Z$ with respect to $\phi,$  
and $Y \xrightarrow{\nu} Z_T$ the normalization, then the smooth locus $U$ of $Y$ over $T$ is dense in all the fibres of $Y$ over $T.$ Note that since the normalization morphism is finite and surjective, $Y$ is a $T$-scheme of relative dimension $d.$ Let $Z_T \xrightarrow{\beta} Z$ be the projection morphism. For each $s\in S,$ we have morphisms between the reduced schemes associated to the fibres $U_{\phi^{-1}(s)} \xhookrightarrow{j} (Y_{\phi^{-1}(s)})_{red} \xrightarrow{(\nu_s)_{red}} ((Z_T)_{\phi^{-1}(s)})_{red} \xrightarrow{(\beta_s)_{red}} (Z_s)_{red},$ such that $(\nu_s)_{red}$ and $(\beta_s)_{red}$ are finite and surjective. Thus, in particular, the composition $\psi:= (\beta_s)_{red} \circ (\nu_s)_{red} \circ j$ is dominant.

We will use the following claim: Let $\mathcal{Y} \xrightarrow{g} \mathcal{Z}$ be a finite
morphism between schemes of pure dimension $d.$ If $V$ is a dense open subset of $\mathcal{Z},$ then $g^{-1}(V)$ is a dense open subset of $\mathcal{Y}.$ Indeed if the open complement $W$ of the closure $\overline{g^{-1}(V)}$ of $g^{-1}(V)$ in $\mathcal{Y}$ is not empty then it is of dimension $d.$ But it maps to the complement of $V$ in $\mathcal{Z},$ which is of dimension $<d.$ This is not possible since $g$ is finite. Thus $W$ is empty and $g^{-1}(V)$ is a dense open subset of $\mathcal{Y}.$

 By generic smoothness on the target applied to $\psi,$ there is a dense open subset $\mathcal{U}$ of $(Z_s)_{red}$ such that the restriction of $\psi$ to $\mathcal{U}$ is smooth. Moreover, by the claim applied to $ (\beta_s)_{red} \circ (\nu_s)_{red},$ we have that the inverse image $\psi^{-1}(\mathcal{U} \cap (Z_s)_{red}^{sm})$ of the dense open subset $\mathcal{U} \cap (Z_s)_{red}^{sm}$ of the smooth locus $(Z_s)_{red}^{sm}$ of $(Z_s)_{red}$ is dense in $U_{\phi^{-1}(s)}.$  
 Thus the restriction of $\alpha$ to  
 $(Z_s)_{red}^{sm}$ vanishes  
 if and only if its restriction to $U_{\phi^{-1}(s)}$ vanishes. Hence the lemma follows from the smooth case, applied to $U$ and the restriction of $\alpha$ to $U.$      
\end{proof}	

\appendix

\section{Symplectic geometry of the cotangent space} \label{subsec:symplectic geometry}

Let $S$ be a scheme and let $Y$ be a smooth $S$-scheme of relative dimension $n$.
Recall that the cotangent space of $Y/S$ is the $Y$-scheme $T^*_{Y/S} \xrightarrow{p_Y} Y:= V((\Omega^1_{Y/S})^*)=\Spec_Y(Sym_{\mathcal{O}_Y}(\Omega^1_{Y/S})^*)$ and hence that the sheaf of germs of $Y$-sections of $T^*_{Y/S}/Y$ is canonically identified with $\Omega^1_{Y/S},$ see \cite[(1.7.9)]{EGAII}. 
Moreover $T^*_{Y/S}$ is a smooth $Y$-scheme of relative dimension $n$ by \cite[Proposition 17.3.8]{EGAIV}, smooth of relative dimension $2n$ as an $S$-scheme. We use the notation $T^*_Y$ (instead of $T^*_{Y/S}$) when there is no risk of confusion, for example when $S$ is the spectrum of a field.
For $f:X\to Y$ a $S$-morphism of smooth $S$-schemes, the pullback of differentials $\Omega^1_{X/S}\xleftarrow{f^*} f^*\Omega^1_{Y/S}$ (\cite[16.4.3.6]{EGAIV}) gives rise to the $X$-morphism $T^*_{X/S} \xleftarrow{(f)_d} X\times_Y T^*_{Y/S}$ called the cotangent map. It is part of the cotangent diagram of $f:$

$$
    \xymatrix{
        T^*_{X/S} & X\times_Y T^*_{Y/S} \ar[l]_-{(f)_d} \ar[d]^{(f)_\pi}\\
          & T^*_{Y/S} }
$$
where $(f)_\pi$ is the canonical projection. 
\begin{rem} \label{rem: pullback section}
Let $U\subset Y$ be an open subset, it follows directly from  
the definitions that if $s_{\alpha}$ is the section of $T^*_{Y/S}/U$ corresponding to $\alpha \in \Gamma(U,\Omega^1_{Y/S}),$ then $(f)_d \circ (X\times_Ys_{\alpha})$ corresponds to $(f^*)^{ad}\alpha\in \Gamma(f^{-1}(U),\Omega^1_{X/S}),$ where $f_*\Omega^1_{X/S}\xleftarrow{(f^*)^{ad}} \Omega^1_{Y/S}$ is adjoint to $f^*.$ \end{rem}

We have the

\begin{lmm} \label{lmm:cotangent of an immersion}
If $f:X\to Y$ is an immersion (resp. a closed immersion) then $(f)_d$ is smooth and surjective and $(f)_\pi$ is an immersion (resp. a closed immersion). Moreover $(f)_d$ admits a section locally on $X.$ 
\end{lmm}
\begin{proof} The morphism $(f)_d$ is smooth and surjective by \cite[Proposition 17.2.5]{EGAIV}, \cite[Proposition 1.7.11 (iii)]{EGAII}, \cite[Proposition 17.3.8]{EGAIV} and stability under base-change of surjective smooth morphisms. We have that $(f)_\pi$ is an immersion (resp. a closed immersion) by \cite[Proposition 4.3.1 (i)]{EGAI}. Finally, local sections of the locally split "conormal" short exact sequence of \cite[Proposition 17.2.5]{EGAIV} induce local sections of $(f)_d,$ by \cite[Proposition 1.7.11 (i)]{EGAII}.
\end{proof}

\begin{defi} \label{defi: canonical form}
The canonical global $S$-relative $1$-form $\theta_{Y/S}$ on the cotangent space of $Y/S$ is the relative $1$-form corresponding to the section $$T^*_{Y/S} \xrightarrow{\Delta_{T^*_{Y/S}/Y}}  T^*_{Y/S}\times_Y T^*_{Y/S} \xrightarrow{(p_Y)_d} T^*_{T^*_{Y/S}/S}$$ of the cotangent space $T^*_{T^*_{Y/S}/S} \xrightarrow{p_{T^*_{Y/S}}} T^*_{Y/S}$, where $\Delta_{T^*_{Y/S}/Y}$ is the diagonal of $T^*_{Y/S} \xrightarrow{p_Y} Y$. 
\end{defi}
Let $\{y_1,...,y_n\}$ be local \'{e}tale coordinates on $Y.$ In terms of the associated local \'{e}tale coordinates $\{y_1,...,y_n;\xi_1,...,\xi_n\}$ on $T^*_{Y/S}$, where $\{\xi_1,...,\xi_n\}$ are dual to $\{dy_1,...,dy_n\}$, $\theta_{Y/S}=\sum_{i=1}^{i=n}\xi_idy_i$. Note that the formation of the canonical form commutes with base-change $S'\to S$ and is compatible with the cotangent diagram, the latter in the sense that $(f)_\pi^*\theta_{Y/S}=(f)_d^*\theta_{X/S}$.

If $S$ is the spectrum of a field $k,$ then we omit the base $S$ from the notations. Let $Q$ be a smooth $k$-scheme of pure dimension $n$. 

\begin{defi} \label{defi: symplectic form} The nondegenerate global exact $2$-form $\omega_Q:=d\theta_Q$ on $T^*_Q$ is called the symplectic form. 
\end{defi}

\begin{defi} \label{defi: Lagrangian}

A subvariety $X \overset{i}\hookrightarrow T^*_Q$ is said to be a Lagrangian subvariety of $(T^*_Q,\omega_Q)$ if it contains a dense open $U\subset X$ on which the symplectic form vanishes, $(i^*\omega_Q)|U=0$ and if at each of its points $x$ it is of dimension $n=dim_xQ$.

\end{defi}

Let us finally note the

\begin{lmm} \label{lmm:isotropy}
Let $f:X\to Y$ be an immersion of smooth $k$-schemes and let $Z_Y \overset{i}\hookrightarrow T^*_Y$ and $Z_X \overset{j}\hookrightarrow T^*_X$ be reduced subschemes. Suppose that $(f)_\pi^{-1}Z_Y=(f)_d^{-1}Z_X$ and that $(f)_\pi^{-1}Z_Y \xrightarrow{(f)_\pi}Z_Y$ is surjective. Then $\omega_Y$ vanishes on a dense open subset of $Z_Y$ if and only if $\omega_X$ vanishes on a dense open subset of $Z_X$.
\end{lmm}
\begin{proof} Note that by Lemma \ref{lmm:cotangent of an immersion}, $(f)_d|_{Z_X}$ is smooth and surjective and $(f)_\pi|_{Z_Y}$ is an immersion. Since by hypothesis $(f)_\pi|_{Z_Y}$ is surjective, it is a nilimmersion.  
Hence, $Z_Y$ being reduced, an isomorphism. Moreover by \cite[Proposition 17.2.3 (ii)]{EGAIV}, \cite[\S 7 n$^{\rm o}$2 Proposition 4]{BourbakiAIII} and flatness of smooth morphisms, the pullback of forms $(f)_d^*:\Omega^2_{Z_X,z} \to \Omega^2_{(f)_d^{-1}Z_X,\tilde{z}}$ is injective for all $z=(f)_d(\tilde{z})$. Since $(f)_\pi^*\omega_Y=(f)_d^*\omega_X$ as $(f)_\pi^*\theta_Y=(f)_d^*\theta_X,$ and  $((f)_d|_{Z_X})\circ((f)_\pi|_{Z_Y})^{-1}$ as well as $((f)_\pi|_{Z_Y})\circ((f)_d|_{Z_X})^{-1}$ preserve open dense subsets, the lemma follows.\end{proof}

\bibliographystyle{alphanum}

\begin{thebibliography}{99}

\bibitem{AusGold} Auslander, M.; Goldman, O. The Brauer group of a commutative ring. Trans. AMS 97 (1960), p. 367-409. 

\bibitem{BerthI} Berthelot, P. $\mathcal{D}$-modules arithm\'{e}tiques I. Op\'{e}rateurs diff\'{e}rentiels de niveau fini. Ann. sc. de l'ENS 29 no. 2 (1996), p. 185-272. 

\bibitem{BerthIntro} Berthelot, P. Introduction \`{a} la th\'{e}orie arithm\'{e}tique des $\mathcal{D}$-modules,
in Cohomologies p-adiques et applications arithm\'{e}tiques II, Ast\'{e}risque 279, p. 1-80 (2002).

\bibitem{BO} Berthelot, P.; Ogus, A. Notes on crystalline cohomology. Princeton University Press, Princeton, N.J.; University of Tokyo Press, Tokyo, 1978. 

\bibitem{BezrBr} Bezrukavnikov, R.; Braverman, A. Geometric Langlands correspondence for $\mathcal{D}$-modules in prime characteristic: the $GL(n)$ case. 
Pure Appl. Math. Q. 3 no. 1 (2007), p. 153-179. 

\bibitem{BMR} Bezrukavnikov, R.; Mirkovic, I.; Rumynin, D. Localization of modules for a semisimple Lie algebra in prime characteristic (with an Appendix by R. Bezrukavnikov and S. Riche: Computation for $\mathfrak{sl}(3)$) Ann. of Math. (2) 167 (2008), no. 3, p. 945-991.

\bibitem{Bjork1} Bj\"{o}rk, J.-E. Rings of differential operators. North-Holland Mathematical Library, 21. North-Holland Publishing Co., 1979. 

\bibitem{Bjork2} Bj\"{o}rk, J.-E. Analytic $\mathcal{D}$-modules and applications. Mathematics and its Applications, 247. Kluwer Academic Publishers Group, 1993.

\bibitem{Borel} Borel, A. et al. Algebraic $D$-modules. Perspectives in Math. 2, Academic Press, 1987.

\bibitem{BourbakiAIII} Bourbaki, N. Alg\`{e}bre, chapitre III in Alg\`{e}bre, Chapitres 1 \`{a} 3 (1970).

\bibitem{BourbakiACI} Bourbaki, N. Alg\`{e}bre commutative, chapitre I in Alg\`{e}bre commutative, Chapitres 1 \`{a} 4 (1969).

\bibitem{DeligneLNM} Deligne, P. Equations diff\'{e}rentielles \`{a} points singuliers r\'{e}guliers. LNM 163.

\bibitem{Hodge II} Deligne, P. Th\'{e}orie de Hodge, II. Publ. Math. de l'IHES, 40 (1971), p. 5-57. 

\bibitem{EisenbudCA} Eisenbud, D. Commutative algebra. With a view toward algebraic geometry. GTM 150, Springer-Verlag, 1995. 

\bibitem{Gabber} Gabber, O. The integrability of the characteristic variety. 
Amer. J. Math. 103 (1981), no. 3, p. 445-468. 

\bibitem{GabberLevasseur} Gabber, O. Equidimensionalit\'{e} de la vari\'{e}t\'{e} caract\'{e}ristique. Expos\'{e} du 18 juin 1982 au s\'{e}minaire parisien sur les alg\`{e}bres enveloppantes.

\bibitem{EGAI} Grothendieck, A. and Dieudonn\'{e}, J. EGA I. Grundlehren der math. Wissenschaften, vol. 166, Springer-Verlag, 1971. 

\bibitem{EGAII} Grothendieck, A. EGA II (r\'{e}dig\'{e} avec la collaboration de Jean Dieudonn\'{e}) \'{E}tude globale \'{e}l\'{e}mentaire de quelques classes de morphismes. Publ. Math. IHES No. 8, 1961.

\bibitem{EGAIV} Grothendieck, A. EGA IV (r\'{e}dig\'{e} avec la collaboration de Jean Dieudonn\'{e}) \'{E}tude locale des sch\'{e}mas et des morphismes de sch\'{e}mas. Publ. Math. IHES No. 20, 24, 28, 32, 1964-1967.

\bibitem{BrauerI} Grothendieck, A. Le groupe de Brauer I. Alg\`{e}bres d'Azumaya et interpr\'{e}tations diverses in Dix Expos\'{e}s sur la Cohomologie des Sch\'{e}mas, p. 46-66. North-Holland, Amsterdam; Masson, Paris, 1968.

\bibitem{FGA Hilb} Grothendieck, A. 
Techniques de construction et th\'{e}or\`{e}mes d'existence en g\'{e}om\'{e}trie alg\'{e}brique IV : les sch\'{e}mas de Hilbert. S\'{e}minaire Bourbaki 221, 1960-1961.

\bibitem{HartshorneAG} Hartshorne, R. Algebraic geometry. GTM 52. Springer-Verlag, 1977. 

\bibitem{Hochschild} Hochschild, G. Simple algebras with purely inseparable splitting fields of exponent 1. 
Trans. AMS 79 (1955), p. 477-489. 

\bibitem{Hoobler} Hoobler, R. Cohomology of purely inseparable Galois coverings. J. Reine Angew. Math. 266 (1974), p. 183-199.

\bibitem{Illusie dR-W} Illusie, L. Complexe de de\thinspace Rham-Witt et cohomologie cristalline. Ann. sc. de l'ENS, 12 no. 4 (1979), p. 501-661.

\bibitem{Kash} Kashiwara, M. D-modules and microlocal calculus. AMS Transl. Math. Monogr., vol. 217, 2003.

\bibitem{KKS} Kashiwara, M.; Kawai, T.; Sato, M. Microfunctions and pseudo-differential equations. Lecture Notes in Math. 287 (1973), Springer-Verlag, p. 264–529. 

\bibitem{Katz} Katz, N. Nilpotent connections and the monodromy theorem: Applications of a result of Turrittin. Publ. Math. de l'IHES, 39 (1970), p. 175-232. 

\bibitem{Kleiman} Kleiman, S. The Picard scheme. Fundamental algebraic geometry, p. 235-321, AMS Math. Surv. and Monogr. vol. 123, 2005.

\bibitem{Kontsevich} 
Kontsevich, M.
Holonomic $\mathcal{D}$-modules and positive characteristic. 
Jap. Journ. of Math. 4 (2009), no. 1, p. 1-25. 

\bibitem{Milne}
Milne, J. \'{E}tale cohomology. Princeton Mathematical Series, 33. Princeton University Press, 1980.

\bibitem{Ogus-Vologodsky}
Ogus, A.; Vologodsky, V.
Nonabelian Hodge theory in characteristic $p$. 
Publ. Math. de l'IHES, 106 (2007), p. 1-138. 

\bibitem{Stacks}
The Stacks project. https://stacks.math.columbia.edu

\bibitem{vdB}
Van den Bergh, M. 
On involutivity of $p$-support, Int. Math. Res. Not. IMRN {\bf 2015}, no.~15, 6295--6304. 

\end{thebibliography}

\end{document}